\documentclass[11pt]{report}

\usepackage{amsmath}
\usepackage{amssymb}

\title{Markov Transitions and the Propagation of Chaos}
\author{Alexander David Gottlieb}

\begin{document}

\maketitle

\begin{abstract}

The propagation of chaos is a central concept of kinetic theory that
serves to relate the equations of Boltzmann and Vlasov to the dynamics
of many-particle systems.  Propagation of chaos means that molecular
chaos, i.e., the stochastic independence of two random particles
in a many-particle system, persists in time, as the number of
particles tends to infinity.

We establish a necessary and sufficient condition for a family of
general n-particle
Markov processes to propagate chaos.  This condition is expressed in
terms of the Markov transition functions associated to the n-particle
processes, and it amounts to saying that chaos of random initial
states propagates if it propagates for pure initial states.

 Our proof
of this result relies on the weak convergence approach to the study of chaos
due to Sznitman and Tanaka.  We assume that the space in which the
particles live is homeomorphic to a complete and separable
metric space so that we may invoke Prohorov's theorem in our proof.

We also show that, if the particles can be in only finitely many states,
then molecular
chaos implies that the specific entropies in the n-particle
distributions converge to the entropy of the limiting single-particle
distribution.

\end{abstract}

\tableofcontents

\def \e {\varepsilon}
\def \pa {\partial}

\def \bx{{\bf x}}
\def \by{{\bf y}}
\def \bv{{\bf v}}
\def \bw{{\bf w}}
\def \bl{{\bf l}}
\def \bs{{\bf s}}
\def \bt{{\bf t}}
\def \bxp{{\bf x'}}
\def \byp{{\bf y'}}
\def \bvp{{\bf v'}}
\def \bwp{{\bf w'}}
\def \bvs{{\bf v^*}}
\def \bws{{\bf w^*}}

\def \ddt{ \frac{d}{dt}}
\def \ppt{ \frac{\partial}{\partial t}}
\def \oon {\frac{1}{n}}
\def \EE {\mathbb{E}}
\def \PP {\mathbb{P}}
\def \RR {\mathbb{R}}
\def \RRd {\mathbb{R}^d}

\def \B {\mathcal{B}}
\def \F {\mathcal{F}}
\def \G {\mathcal{G}}
\def \I {\mathcal{I}}
\def \J {\mathcal{J}}
\def \P {\mathcal{P}}
\def \T {\mathcal{T}}
\def \Kt {\widetilde{K}}
\def \phat {\widehat{\phi}}

\def \s {\sigma}

\newtheorem{proposition}{Proposition}[chapter]
\newtheorem{theorem}{Theorem}[chapter]
\newtheorem{lemma}{Lemma}[chapter]
\newtheorem{definition}{Definition}[chapter]
\newtheorem{corollary}{Corollary}[chapter]

\chapter{Introduction}
\label{Introduction}
\section{Overview}

  Kinetic theory is the analysis of
nonequilibrium physical phenomena that emerge from the collective
behavior of large numbers of particles.   That analysis is
accomplished by the techniques of probability theory; kinetic
theory has an inherently statistical character.  One of the
notions of probability theory from which one can derive Boltzmann's
equation and Vlasov's equation, two staples of kinetic theory, is the
propagation of chaos.

The concept of propagation of chaos originated with Kac's Markovian models
of gas dynamics \cite{Kac}.   Kac invented a class of interacting particle
systems wherein particles collide at random with
each other while the {\it density} of particles evolves
deterministically in the limit of infinite particle number.
A nonlinear evolution equation analogous to Boltzmann's
equation governs the particle density. Gr\"unbaum proved the
propagation of chaos along Kac's lines for the spatially homogeneous
 Boltzmann equation given existence
and smoothness assumptions on the Boltzmann semigroup \cite{Gru}.
The processes of Kac were further investigated with regard to
their fluctuations about the deterministic infinite particle limit
in \cite{McK75, Uch83, Uch88}.

 McKean introduced propagation
of chaos for interacting diffusions and analyzed what are now called
McKean-Vlasov equations \cite{McK66, McK68}.  Independently,
Braun and Hepp \cite{BH} analyzed the propagation of chaos for Vlasov
equations and proved a central limit theorem for the fluctuations.
Analysis of the fluctuations and large deviations for McKean-Vlasov
processes was carried out in \cite{Tanaka, Sz84, DawsonGartner}.
Chorin~\cite{Chorin} created a numerical method for the two
dimensional Navier-Stokes equation by interleaving independent random
walks into the discretized dynamics of interacting {\it vortex blobs},
smoothed and localized patches of vorticity that move without changing shape.
Propagation of chaos has been studied in connection with this {\it
random vortex method} by  \cite{MarPul, Osada, Meleard2}.
Other instances of the propagation of chaos have been studied in
\cite{Oelsch, Rez, RezTar, Graham}.   A thorough analysis of the
convergence of numerical schemes based on stochastic particle methods
for McKean-Vlasov equations in one dimension is undertaken in
\cite{BossyTalay, Talay}.

Finally, we refer the reader to the long, informative articles by
 Sznitman~\cite{Sznitman} and by
M\'el\'eard~\cite{Meleard} in Springer-Verlag's {\it Lecture Notes in
Mathematics}.

The aforementioned authors are mostly concerned with
proving that specific systems propagate chaos, rather than the
propagation of chaos {\it per se}.  The modest purpose of this
dissertation is to clarify the definition of propagation of
chaos in general, and not to prove that any particular system propagates chaos.
The essential content of this dissertation is
Definition~\ref{PoCDef1} and Theorem~\ref{theorem2} of
Chapter~\ref{Propagation of Chaos}.

This dissertation is organized as follows.

  The rest of Chapter \ref{Introduction} is an informal
summary of our point of view and contains a statement of our main theorem.
We introduce general Markovian interacting particle systems and adopt a
strong-sense definition of the propagation of chaos.  We can then
characterize the propagation of chaos in terms of
the Markov transition functions that define the interacting particle
systems.

Chapter \ref{Kinetic Theory} describes the most important instances of the
propagation of chaos.
The concept of propagation of chaos is most useful in (and was indeed
motivated by) the kinetic theories of gases, plasmas, and stellar systems.
Boltzmann's equation for dilute
gases is discussed in Section~\ref{Dilute Gases} and Vlasov's equation
for plasmas and stellar systems is discussed in Section~\ref{Plasmas}.

Chapters \ref{Chaos and Weak Convergence} and \ref{Propagation of Chaos}
are meant to be self-contained, formal, and brisk.  They contain the
necessary background
and the proofs of the theorems that flow from the point of view
described in Section~\ref{Definition of Propagation of Chaos and Statement
of Main
Result}.

Chapter~\ref{Chaos and Weak Convergence} discusses the
theorem of Sznitman and Tanaka, which is
our main technical tool.  A detailed proof of this theorem is given in
Section~\ref{Theorem of Sznitman and Tanaka}.
Theorem~\ref{SpecEntThm} states that, if the underlying space is
finite, $p$-chaos entails the convergence of specific entropy to the
entropy of $p$.

Chapter~\ref{Propagation of Chaos} is dedicated to the proof of our
main theorem and its corollaries.  Although the theorems there are neither deep
nor surprising, they should still be of interest to those who work with
the propagation of chaos because they establish
properties of the propagation of chaos which, though easy to take for
granted, do require some proof.  The proofs rely
on the theorem of Sznitman and Tanaka and basic properties of
convergence in law.  Our main theorem requires the completeness of the
basic space so that we may invoke Prohorov's theorem in its proof.

This work was supported by the U.S. Department of Energy at the
Lawrence Berkeley National Laboratory, at the Mathematics Department,
under Alexandre Chorin.

\section{Definition of Propagation of Chaos and Statement of Main Result}
\label{Definition of Propagation of Chaos and Statement of Main Result}

Statistical mechanics and kinetic theory are probabilistic theories of
many-body systems;
their predictions are intended to be valid only when the number of
particles is very large, typically as large as Avogadro's number.
The equations of kinetic theory are obtained by studying the limiting
behavior of $n$-particle systems as $n$ tends to infinity.  A key
concept in such studies is the propagation of chaos.

 The concept of propagation of chaos was motivated originally by the kinetic
theories of gases and plasmas.
  Before we delve into these kinetic theories, we first set
up, in this section, a very general framework for the
study of interacting particle systems and the propagation of chaos.

Our particles shall live in a space $S$.  For gases and plasmas $S$
would be position-velocity space, a subset of $\RR^6$, but for our general
purposes $S$ may be any
separable metric space.  The state  of an
$n$-particle system is a point in $S^n$, the $n$-fold Cartesian
product of $S$ with itself, also a metric space.  Whether the dynamics of
the $n$-particle system are deterministic or random, we desire that
the future motion of a system of particles depend only on its current
state, and not the entire history of the particles' motion.  We stipulate that
the $n$-particle dynamics are Markovian; the future depends on the
past only through the present state.  Markovian processes are defined
by their transition functions, so our data includes one transition
function $K_n(\bs,B, t)$ for each $n$.  That is, for each $n$, we
are given the transition function
\[
   K_n: S^n \times \B_{S^n} \times [0,\infty) \longrightarrow [0,1] ,
\]
where $\B_{S^n}$ is the Borel $\s$-algebra on $S^n$.  The transition functions
$K_n$ have the following interpretation.  For $t \ge 0$,
$\bs \in S^n$, and $B \in \B_{S^n}$, the probability that the state at time
$t$ of
an $n$-particle system belongs to $B$, given that the state was
initially $\bs$, is $K_n(\bs,B,t)$.   The Markov property implies that
the transition functions satisfy the Chapman-Kolmogorov equations:
\[
      K_n(\bs,B,t+t') = \int_{S^n} K(\bs,d{\bf s'},t') K({\bf s'},B,t)
\]
for all $ t,t' \ge 0, \bs \in S^n,$ and $ B \in \B_{S^n}$.

We restrict our consideration to $n$-particle systems whose dynamics
commute with permutations by imposing the following conditions on
the transition functions.   Let $\Pi_n$ denote the set of permutations of
$\{ 1,2,\ldots,n \}$.  If $\pi \in \Pi_n$ and $\bs = (s_1,s_2,\ldots,s_n)
\in S^n$, let
\[
    \pi \cdot \bs = (s_{\pi(1)},s_{\pi(2)},\ldots,s_{\pi(n)}) ,
\]
and, for $B \subset S^n$, let
\[
     \pi \cdot  B = \left\{ \pi \cdot \bs: \bs \in B \right\}.
\]
We suppose that the transition functions satisfy
\begin{equation}
     K_n(\bs,B,t) = K_n(\pi \cdot \bs, \pi \cdot  B, t)
\label{permutations}
\end{equation}
for all permutations $\pi$, points $\bs$, Borel sets $B$, and times $t$.

Here, then, is the set-up.   A separable metric space $S$ is given,
along with a sequence
\[
       \left\{ K_n(\bs,B,t)  \right\}_{n=1}^{\infty}
\]
of Markov transition functions that satisfy the permutation
condition (\ref{permutations}), $K_n$ being a transition function on
$S^n$.  We have a Markovian dynamics of $n$
particles, an $n$-particle system, for each $n$.
 Propagation of chaos, defined shortly, is an attribute of {\it
 families} of particle systems, indexed by $n$;
a family of $n$-particle systems either does or does not propagate chaos.

In order to define the propagation of chaos we must first define the
property of being chaotic --- ``chaos'' for short.  For
each $n$, let $\rho_n$ be a symmetric probability measure on $S^n$,
i.e., a probability measure on $S^n$ such that
\[
      \rho_n(\pi \cdot  B) = \rho_n(B)
\]
for all permutations $\pi$ and all $B \in \B_{S^n}$.  Let $\rho$ be a
probability measure on $S$.

\smallskip
\noindent {\bf Definition}:
{\it
The sequence $\{\rho_n \}$ is $\rho$-{\bf chaotic} if, for any natural
number $k$ and any bounded continuous functions
$
    g_1(s), g_2(s),\ldots,g_k(s)
$
on $S$,
\[
    \lim_{n \rightarrow \infty} \int_{S^n} g_1(s_1) g_2(s_2) \cdots
    g_k(s_k) \rho_n(ds_1 ds_2 \cdots ds_n) = \prod_{i=1}^k \int_S
    g_i(s) \rho(ds)  .
\]
}
\smallskip
In words, a sequence probability measures on the product
spaces $S^n$ is $\rho$-chaotic if, for fixed $k$, the joint
probability measures for the first $k$ coordinates tend to the product
measure $\rho(ds_1)\rho(ds_2)\cdots\rho(ds_k) \equiv \rho^{\otimes k}$
on $S^k$.  If the measures $\rho_n$
are thought of as giving the joint distribution of $n$ particles
residing in the space $S$, then $\{ \rho_n \}$ is $\rho$-chaotic if
$k$ particles out of $n$ become more and more independent as
$n$ tends to infinity, and each particle's distribution tends to
$\rho$.  A sequence of symmetric probability measures on $S^n$ is {\it
chaotic}
if it is $\rho$-chaotic for some probability measure $\rho$ on $S$.

If a Markov process on $S^n$ begins in a random state with
distribution $\rho_n$, the distribution of the  state after $t$ seconds
of Markovian random motion
can be expressed in terms of the transition function $K_n$ for the
Markov process.   The distribution at time $t$ is the probability
measure $ U^n_t \rho_n $ defined by
\begin{equation}
\label{U}
     U^n_t \rho_n (B) := \int_{S^n} K_n(\bs, B, t) \rho_n(d\bs)
\end{equation}
for all $B \in \B_{S^n}$.  If $K_n$ satisfies the permutation
condition (\ref{permutations}) then $ U^n_t \rho_n $ is symmetric
whenever $\rho_n$ is.

\smallskip
\noindent {\bf Definition}:
{\it
A sequence
\[
       \left\{ K_n(\bs,B,t) \right\}_{n=1}^{\infty}
\]
whose $n^{th}$ term is a Markov transition function on $S^n$ that
satisfies the permutation condition
(\ref{permutations}) {\bf propagates chaos} if, whenever
$\{\rho_n\}$ is chaotic so is $ \{ U^n_t \rho_n \} $
for any $t \ge 0$, where
$U^n_t$ is as defined in (\ref{U}).
}
\smallskip

We sometimes say that a family of $n$-particle Markov {\it processes}
propagates chaos when we really mean that the associated family of
transition functions
propagates chaos.

  It follows
from the definition of propagation of chaos that for each $t>0$ there
exists an operator $ U_t^{\infty}$
on probability measures such that
$\{ U^n_t \rho_n \}$ is $ U_t^{\infty}\rho$-chaotic if $\{ \rho_n \}$
is $\rho$-chaotic.  This operator is typically nonlinear; even though
$U_t^{\infty}$ is derived from the linear operators $U_t^n$ by taking
a limit of sorts, it is not actually a limit of linear operators, and may
be nonlinear.
For families of interacting particle systems suited
to the study of gases or plasmas, the semigroup
$\{ U_t^{\infty} \}_{t \ge 0} $ is
 the semigroup of solution operators for the Boltzmann or the Vlasov
equation.  (The existence of the operators $ U_t^{\infty} $
is part of our main theorem, stated shortly.)

We are adopting here a strong definition of the propagation of chaos.
Other authors \cite[p. 42]{Meleard}\cite[p. 98]{Pulvirenti}
 have defined propagation of chaos in a weaker
sense:  a family of Markovian $n$-particle processes propagates chaos
if  $\{U_t^n \rho^{\otimes n}\}$ is chaotic for all
$\rho \in \P(S)$  and $t > 0$, where $ \rho^{\otimes n}$ is
product measure on $S^n$.   For these authors,
only {\it purely} chaotic sequences of initial measures
are required to ``propagate'' to chaotic sequences.   This condition
is strictly weaker than the one we adopt for our definition.
For example, take $S = \{0,1\}$ and let $\delta(x)$ or $\delta_x$
denote a point mass at $x$.  Then, if
\[
   K_n(\bs,\cdot,t) = \left\{ \begin{array}{cc}
                 \delta_{(1,1,\ldots,1)} & \mathrm{if}\  \bs \ne
(0,0,\ldots,0) \\
                 \delta_{(0,0,\ldots,0)} & \mathrm{if}\  \bs =
(0,0,\ldots,0)   \\
                  \end{array} \right.
\]
for all $t>0$, the sequence $\{K_n\}$ propagates chaos in the weak sense, but
not in the strong sense of our definition.   Under these $K_n$'s,
 the $\delta(0)$-chaotic sequence $\{ \delta_{(0,0,\ldots,0)} \}$ is
propagated to itself, while other $\delta(0)$-chaotic sequences
are propagated to $\delta(1)$-chaotic sequences, and yet other
$\delta(0)$-chaotic sequences are not propagated to chaotic sequences
at all.

Our main result is a
condition on the Markov transition functions for a family of
$n$-particle processes that is necessary and sufficient for the
propagation of chaos (in the strong sense).
Before we state it we must recall the weak topology on probability
measures and introduce some necessary notation.

If $X$ is a completely regular topological space (as normal
topological spaces are), let $\P(X)$ denote
the space of probability measures on $X$ endowed with the weakest topology
relative to which all the functions $I_g:\P(X) \longrightarrow \RR$
are continuous, where
\[
    I_g(\mu) = \int_X g(x)\mu(ds)
\]
and $g$ ranges over the bounded and continuous real-valued functions on $X$. A
sequence $\{\mu_n\}$ in $\P(X)$ converges to $\mu$ in this weak topology if
\[
     \int_X g(x)\mu_n(ds) \longrightarrow \int_X g(x)\mu(dx)
\]
for all $g \in C_b(X)$, the space of bounded and continuous real-valued
functions on $X$.

For $\nu$ a measure on $S^n$, let $\widetilde{\nu}$ denote its
symmetrization:  for all $B \in \B_{S^n}$,
\[
    \widetilde{\nu}(B) = \frac{1}{n!}\sum_{\pi \in \Pi_n} \nu(\pi \cdot  B) .
\]
For fixed $\bs \in S^n$ and $t \ge 0$,
denote by $\Kt_n(\bs, \cdot, t)$ the symmetrization of the
measure $K_n(\bs, \cdot, t)$.

For
\[
  \bs = (s_1,s_2,\ldots,s_n) \in S^n,
\]
let $\e_n(\bs)$ denote the purely atomic probability measure
\begin{equation}
    \e_n( \bs) = \oon \sum_{i=1}^n \delta(s_i) \quad ,
\label{e}
\end{equation}
where $\delta(s)$ --- Dirac's delta --- denotes a point-mass at $s$.
The function $\e_n$ takes ordered $n$-tuples to purely atomic probability
measures
consisting of $n$ point-masses of weight $\oon$ each.

We can now state our main theorem.

\smallskip
\noindent {\bf Main Theorem}:
{\it
Let $S$ be a complete, separable metric space.
Let the Markov transitions $K_n$ satisfy the permutation condition
(\ref{permutations}).

  Propagation of chaos by the sequence $\{K_n\}$ is equivalent to the
following condition:

For every $t>0$ there exists a continuous map
\[
     U_t^{\infty}: \P(S) \longrightarrow \P(S)
\]
such that, if the sequence
\[
    \bs_1 \in S, \quad \bs_2 \in S^2, \quad \bs_3 \in S^3, \ldots
\]
is such that $\{ \e_n(\bs_n) \}_{n=1}^{\infty}$ converges to $\rho$
in $\P(S)$, then the sequence of symmetric measures
\[
     \left\{ \Kt_n(\bs_n,\cdot,t) \right\}_{n=1}^{\infty}
\]
is $U_t^{\infty}\rho$-chaotic.
}
\smallskip

The necessity of the condition of the preceding theorem is an easy
consequence of
the definition; its sufficiency is nontrivial.  Our theorem
shows that to prove propagation of chaos it is {\it sufficient} to verify that
$\{ U^n_t\rho_n \}$ is chaotic when the initial measures $\rho_n \in
\P(S^n)$ are symmetric atomic measures of the form
\begin{equation}
\label{specialmeasures}
       \rho_n = \frac{1}{n!}\sum_{\pi \in \Pi_n}
       \delta(\pi\cdot\bs_n); \quad \bs_n \in S^n.
\end{equation}
This sufficient condition can come in handy.
In 1977, Braun and Hepp \cite{BH} proved the propagation of chaos for Vlasov's
equation, provided the initial conditions are ``pure initial states'' of
the form (\ref{specialmeasures}).
Sznitman \cite{Sz84}, in 1983, noted that Braun and Hepp
require ``purely atomic initial data'' to propagate their
chaos, implicitly suggesting that this restriction to special initial
conditions weakens their result.  Our theorem
shows that it did indeed suffice for Braun and Hepp to verify propagation
of chaos for purely
atomic initial data.

The theorem is proved by expressing chaos in terms of weak convergence of
probability measures in $\P(\P(S))$ and then applying Prohorov's theorem.
 Prohorov's theorem \cite{Billingsley} states that a family $\F$ of
probability measures on a complete and separable metric space is
relatively compact if {\it and only if} $\F$ is tight.
Our hypothesis that $S$ is complete ensures that $\P(S)$ is also
complete and enables us to apply Prohorov's theorem in $\P(\P(S))$.
Chapter \ref{Propagation of Chaos} is devoted to the proof of the theorem.

The study of chaos via weak convergence in the space
$\P(\P(S))$ is due to Sznitman
\cite{Szn} and Tanaka \cite{Tanaka}.  They proved that a sequence of
symmetric measures
$\{ \rho_n \}$ is $\rho$-chaotic if and only if the
probability measures induced on $\P(S)$ by
$ \e_n $ converge in $\P(\P(S))$ to $\delta(\rho)$,  a point mass at $\rho
\in \P(S)$.
This device is essential to our approach and is discussed at length in
Chapter \ref{Chaos and
Weak Convergence}.

\chapter{Kinetic Theory and the Propagation of Chaos}
\label{Kinetic Theory}

Boltzmann's equation for dilute gases and Vlasov's equation for
plasmas govern the evolution, the change over time, of the density of
particles in position-momentum space.   The particle density changes
due to interactions between the particles: binary collisions of
molecules in a dilute gas or mutual electric forces acting between ions in a
plasma.  The rate of change of the particle density
is determined by the particle density {\it itself} through the
particle interactions.   The evolution equations
of Boltzmann and Vlasov are nonlinear because of the way the particle
density affects its own evolution.

  This chapter reviews the equations
of Boltzmann and Vlasov for the sake of illuminating the meaning and
physical relevance of the propagation of chaos.
One may consult \cite{Spohn} for a more thorough treatment of kinetic
theory.

Section~\ref{Dilute Gases} presents the theory of dilute gases from the
point of view of the propagation of chaos.  First, the classic derivation of
Boltzmann's equation is repeated in \ref{Boltzmann's equation}.  Then,
in \ref{Particle systems for Boltzmann's equation}, two types of
$n$-particle systems are introduced that satisfy Boltzmann's equation
in the infinite particle limit.

Section~\ref{Plasmas} is about Vlasov's equation for plasmas and stellar
systems.  Vlasov's equation is introduced in \ref{Vlasov's equation}
and rederived in terms of the propagation of chaos in
\ref{McKean-Vlasov particle systems}.

\section{Dilute Gases}
\label{Dilute Gases}

\subsection{Boltzmann's equation}
\label{Boltzmann's equation}

In this section we summarize Boltzmann's derivation of his equation
for a dilute gas.  Our source is the first chapter of his {\it Lectures on
Gas Theory}
\cite{Boltzmann}, written over a century ago.

Boltzmann modeled the
molecules of the gas by hard spheres: balls of radius $r$ that collide
elastically according to simple mechanics.  When a ball having
velocity $\bv$ collides with a ball having velocity $\bw$, the
collision instantaneously changes the velocity of the first ball
from $\bv$ to $\bvp$ and the velocity of the second ball from $\bw$
to $\bwp$.   Given the relative orientation of the balls at the time
of impact,
the, post-collisional or outgoing velocities are determined by
the laws of conservation of energy and momentum.
  Suppose that, at the moment of impact, $\bl$ is the unit vector
parallel to the ray that originates at the center of the ball of
velocity $\bv$ and passes through the center of the ball of
velocity $\bw$.  Such a collision, which we call a $(\bv,\bw:
\bl)$ collision, changes the velocities of the balls to
\begin{eqnarray}
\label{collisions}
   \bv &\longrightarrow & \bvp = \bv + [(\bw - \bv)\cdot \bl]\bl  \nonumber
\\
   \bw &\longrightarrow & \bwp = \bw - [(\bw - \bv)\cdot \bl]\bl    .
\nonumber \\
\end{eqnarray}
A collision of type $(\bv,\bw: \bl)$ is only possible if $(\bw - \bv)\cdot
\bl <
0$.  Except during collisions, which have instantaneous duration,
molecules (hard spheres in this model) travel
inertially, with unchanging velocity.
 Let $n$ denote the number of molecules in the
gas, and let the number of molecules per unit volume of position-momentum
space be
given by the density
\[
         f(\bx,\bv,t)d\bx d\bv ,
\]
so that the proportion of molecules which, at time $t$, are located in
a region $X$ of space and have velocities belonging to a set $V$ of
velocities is
\[
   \frac{1}{n} \int_V \int_X f(\bx,\bv,t)d\bx d\bv  .
\]
Boltzmann's equation tells how $f(\bx,\bv,t)$ changes due to the
collisions detailed above.

 The density $f(\bx,\bv,t)d\bx d\bv$ of molecules
changes through the inertial motion of the molecules between
collisions (called {\it free streaming}) and through collisions between
molecules.  Boltzmann's equation can be written
\[
   \ppt f(\bx,\bv,t) = - \bv \cdot \nabla_{\bx}f(\bx,\bv,t) +
   Q[f(\bx,\bv,t)] ,
\]
where $- \bv \cdot \nabla_{\bx}f$ gives the rate of change of $f$ due
to free streaming, and $Q[f]$, the collision operator applied to $f$, gives
the rate of
change of the density due to collisions.

Further assumptions are needed to determine $Q[f]$, the rate of change of
$f(\bx,\bv,t)$ due to collisions.
  We know the effect of a $(\bv,\bw: \bl)$ collision, but we
also need to know the rate at which those collisions are occurring.
Boltzmann assumed that the rate at which $(\bv,\bw: \bl)$ collisions are
happening
at a point $\bx$ of space is proportional to
$r^2 \| (\bw - \bv)\cdot \bl \|$ and jointly proportional to
the densities at $\bx$ of molecules having velocities $\bv$ and $\bw$.
These assumptions are the {\it Stosszahlansatz},
or collision-number-hypothesis:  the rate of $(\bv,\bw: \bl)$ collisions
at $\bx$ is
\begin{equation}
\label{Stosszahlansatz}
   r^2 \| (\bw - \bv)\cdot \bl \|f(\bx,\bv,t) f(\bx,\bw,t).
\end{equation}

The rate of change of $f(\bx,\bv,t)$ due to collisions, $Q[f]$, equals
the rate at which the molecules are receiving post-collisional velocities
$\bv$ less the rate at which molecules already having velocity $\bv$
are colliding with other molecules and exchanging $\bv$ for other
velocities.  The loss rate is easy to express,
assuming the Stosszahlansatz:
\begin{equation}
   L[f] := \frac {r^2}{2} \int_{\RR^3} \int_{S_2} f(\bx,\bv,t)
   f(\bx,\bw,t)  \| (\bw - \bv)\cdot \bl \| d\bl d\bw   ,
\label{loss}
\end{equation}
where $S_2$ is the unit sphere in $\RR^3$ and $d\bl$
indicates the normalized and uniform measure on the sphere $S_2$, is the
 number of molecules per unit
volume at $\bx$ of velocity $\bv$
  that will collide with other molecules between times
$t$ and $t + \Delta t$, divided by $\Delta t$.

There is a similar expression for the gain rate at which collisions are
{\it resulting} in molecules having velocity $\bv$.  Observe that a binary
collision can only
produce a post-collisional, or outgoing, velocity $\bv$ if the
velocities before collision were $\bv + ((\bw - \bv)\cdot \bl)\bl$ and
$\bw - ((\bw - \bv)\cdot \bl)\bl$ for some $\bw$.   Let
\begin{eqnarray*}
     \bvs &=& \bv + ((\bw - \bv)\cdot \bl)\bl    \\
     \bws &=& \bw - ((\bw - \bv)\cdot \bl)\bl    .\\
\end{eqnarray*}
The number of
molecules per unit volume that will end up having velocity $\bv$
because of a collision that took place between times $t$ and
$t + \Delta t$, divided by $\Delta t$, equals
\begin{equation}
   G[f] := \frac{r^2}{2} \int_{\RR^3} \int_{S_2} f(\bx,\bvs,t)
   f(\bx,\bws,t) \| (\bw - \bv)\cdot
   \bl \| d\bl d\bw .
\label{gain}
\end{equation}
In fact, $\bvs = \bvp$ and $\bws =\bwp$; if a $(\bv, \bw:\bl)$
collision changes $\bv$ to $\bvp$ and $\bw$ to $\bwp$, then a $(\bvp,
\bwp:\bl)$
collision changes $\bvp$ to $\bv$ and $\bwp$ to $\bw$.
It is only a lucky accident that $\bvs = \bvp$, so we
 emphasize, by introducing new notation, that $\bvs$ and $\bws$ are
supposed to be
velocities for which a $(\bvs, \bws:\bl)$ collision {\it results in} a
velocity $\bv$.

The net rate of change of $f(\bx,\bv,t)$ due to collisions equals the
gain rate minus the loss rate:  $Q[f] = G[f] - L[f]$.  Boltzmann's
equation is thus
\begin{equation}
      \ppt f(\bx,\bv,t) + \bv \cdot \nabla_{\bx}f(\bx,\bv,t) =
   G[f(\bx,\bv,t)] - L[f(\bx,\bv,t)],
\label{BoltzEq}
\end{equation}
where $G[f]$ and $L[f]$ are as defined in (\ref{gain}) and
(\ref{loss}).

The existence of solutions of Boltzmann's equation
(\ref{BoltzEq}) is difficult to prove.  The state of the art is the
global existence of mild solutions proved by Di Perna and Lions \cite{DipLio}.

\subsection{Particle systems for Boltzmann's equation}
\label{Particle systems for Boltzmann's equation}

Kac~\cite{Kac55},
in his article {\it Foundations of Kinetic Theory} of 1954,
propounds the relationship between Boltzmann's equation and certain
$n$-particle Markovian jump processes.
These $n$-particle systems are inherently stochastic; the collisions
have random results and happen at random times.  The dynamics are not
the true dynamics of deterministically colliding molecules, rather, the
stochastic motion of fictitious particles which
obey the spatially homogeneous Boltzmann equation on the
macroscopic level.

The spatially homogeneous Boltzmann equation is the equation satisfied by
a position-velocity density that does not depend on position:
$f(v)dv$.  So Kac imagines a gas of $n$ particles on the line,
particles whose positions are unimportant and are not given, but whose
velocities
\begin{equation}
\label{velocities}
     v_1,\  v_2, \ldots,\  v_n; \qquad v_i \in \RR
\end{equation}
completely specify the state of the gas.  Kac proposes a stochastic
dynamics of
these states driven by collisions between pairs of particles.  Suppose
that the state is initially given by the list (\ref{velocities}).
At a random time, a collision occurs.  A collision changes the values
of a random pair of the $n$ velocities in the list, at random.  Once
the state of the gas has jumped to a new state due to a collision,
another random time elapses, another collision occurs, and so forth.
The random times are taken to be independent and to have exponential
distributions with mean duration $\tau /n$; the probability that a
collision happens later than $t$ seconds after the previous
collision is $e^{- n t / \tau}$.  Notice that the more particles
there are, the faster collisions are occurring.  Each collision only
affects the
velocities of two particles, the affected pair
being selected at random from one of the $n(n-1)/2$ possible pairs of
particles.
Given that a pair of particles having velocities $v$ and $w$ collide,
those two velocities change to another pair $v'$ and $w'$
satisfying the conservation of energy condition
\[
       (v')^2 + (w')^2 = v^2 + w^2 ,
\]
but otherwise at random, so that $(v',w')$ is randomly sampled from the
uniform probability measure on
the circle
\[
      \left\{ (v',w'): (v')^2 + (w')^2 = v^2 + w^2 \right\}  .
\]
Kac's $n$ particle gas is thus a Markov jump process
on $\RR^n$, for each $n$.

In \cite{Kac, Kac55}, Kac proves that this family of
$n$-particle gases propagates chaos.  Indeed, the exact definition of
chaos as the asymptotic independence of particles is due to Kac.  The
notion of chaos originates in Boltzmann \cite{Boltzmann}, who derived his
equation under a hypothesis of ``molecular disorder (chaos).''

Kac proved that if the particles of
each $n$-particle gas initially have independent and
$f_0(v)dv$-distributed velocities, then at a later time $t$ the velocities
of a random pair become increasingly independent as $n \longrightarrow
\infty$, even though the initial
condition of {\it pure} independence or ``molecular chaos'' has been
spoiled by
collisions.   The random velocity of a single particle at time
$t$ becomes increasingly $f(v,t)$-distributed as $n \longrightarrow
\infty$, where $f(v,t)$ satisfies an analog of Boltzmann's equation, namely,
\begin{eqnarray*}
     \ppt f(v,t) &=& \frac{2}{\tau} \int_{\RR} \int_0^{2\pi} f( v \cos
     \theta - w \sin\theta) f( v \sin\theta + w \cos\theta) d\theta dw
     -f(v) \\
     f(v,0) &=& f_0(v)    . \\
\end{eqnarray*}
Indeed, the sequence of
$n$-particle joint distributions at time $t$ is $f(v,t)dv$-chaotic.

Similar procedures yield particle systems for the
spatially homogeneous Boltzmann equation \cite{Kac55}.
The spatially homogeneous Boltzmann equation for hard spheres of
radius $r$ is
\begin{eqnarray}
   \ppt f(\bv,t)  &=&
        G[f(\bv,t)] - L[f(\bv,t)]     \nonumber  \\
   G[f] &=& \frac{r^2}{2} \int_{\RR^3} \int_{S_2} f(\bvs,t)
     f(\bws,t) \| (\bw - \bv)\cdot
     \bl \| d\bl d\bw   \nonumber  \\
   L[f] &=& \frac {r^2}{2} \int_{\RR^3} \int_{S_2} f(\bv,t)
     f(\bw,t)  \| (\bw - \bv)\cdot \bl \| d\bl d\bw   , \nonumber \\
\label{SHBE}
\end{eqnarray}
where $d\bl$ is normalized surface area on the sphere $S_2$, and
\begin{eqnarray}
     \bvs &=& \bv + ((\bw - \bv)\cdot \bl)\bl  \nonumber  \\
     \bws &=& \bw - ((\bw - \bv)\cdot \bl)\bl  \nonumber  .\\
\label{stars}
\end{eqnarray}
One may devise several $n$-particle jump processes for the Boltzmann
equation.
 Gr\"unbaum \cite{Gru} suggests
one with a three-stage random mechanism for making jumps:  given that the
initial
state of the gas or the state it has just jumped to is
\[
      (\bv_1, \bv_2, \ldots, \bv_n); \qquad \bv_i \in \RR^3,
\]
\begin{trivlist}
\item{1)} Select two distinct particles at random (equiprobably), say
the $i^{th}$ and $j^{th}$ particles where $i < j$.
\item{2)} If $\bv_i = \bv_j$ select another pair.  Otherwise
wait for an exponentially distributed random time of mean duration $\|
\bv_i - \bv_j \| / (n-1) $.
\item{3)}  Jump  to $(\bv_1, \ldots, \bvs_i,\ldots, \bvs_j,  \ldots,
\bv_n)$
with probability proportional to
$\frac {\| (\bv_i - \bv_j)\cdot \bl \| } {\| \bv_i - \bv_j \| }$, where
$\bvs_i,\bvs_j$, and $\bl$
are as in (\ref{stars}).
\end{trivlist}
Note that the jumps speed up as the number of
particles increases so that the number of jumps per particle per unit time
is roughly constant.  Gr\"unbaum proves that this family of
$n$-particle processes propagates chaos and that the limit satisfies
(\ref{SHBE}) under certain assumptions \cite{Gru}.  His proof relies
on the theory of strongly continuous contraction semigroups.

Other jump processes similar to those of Kac have been treated by several
authors.
Uchiyama \cite{Uch88} proves propagation of chaos and a central
limit theorem for families of Kac-type processes, on {\it countable} sets
of velocities.
Rezakhanlou and Tarver \cite{RezTar} prove an interesting propagation
of chaos result for the discrete Boltzmann equation in one dimension.
Their particles travel with constant velocities around a circle in
between random collisions that become increasingly local as the number
of particles increases.
Graham and M\'el\'eard \cite{Meleard} prove the propagation of chaos
for a variant of the Boltzmann equation with nonlocal collisions.  Particles
experience random Kac-type collisions, but do not need to be at the
same spatial location in order to collide.  Bird's numerical scheme
for Boltzmann's equation \cite{Bird} amounts to the simulation of
one of the processes studied by Graham and M\'el\'eard.

The jump processes of Kac et alia are intrinsically stochastic, for
collisions happen at random and have random results.  On the other
hand, the dynamics of real molecules are strictly deterministic, or
are classically conceived to be such.  Our idealized model for
molecular dynamics, the hard sphere model, admits no randomness at
all.  When two particles collide, their outgoing velocities are
determined by their incoming velocities and their attitude at
collision.  (It is true that the outcome of a simultaneous collision of
three or more spheres may not be determined, but in a dilute enough
hard sphere gas these collisions are so rare that they have negligible
effect.)  Since Boltzmann's equation is supposed to govern the
macroscopic behavior of the density of a hard sphere gas, it ought to
be derivable somehow from the deterministic dynamics of hard spheres.
But alas, it would appear that the Boltzmann equation is not even {\it
consistent} with molecular dynamics, much less derivable from it, for
the molecular dynamics are reversible and Boltzmann's equation is
irreversible.  This apparent antinomy, known as Loschmidt's paradox,
has been raising deep concerns about the validity of Boltzmann's
equation for nearly as long as that equation has been known.  It is
therefore surprising and philosophically significant that
(notwithstanding Loschmidt's paradox)
Boltzmann dynamics can indeed be derived from molecular dynamics.

  Grad \cite{Grad} first advanced the idea that Boltzmann's equation may be
derived in the {\it dilute limit}
\[
    nr^2\longrightarrow \mathrm{constant}
\]
of hard sphere dynamics, and
 Lanford \cite{Lanford} succeeded in a rigorous derivation of Boltzmann's
equation
along the lines suggested by Grad.
 Lanford's theorem can be neatly expressed in terms of of chaos.
This approach can be found in {\it The Mathematical Theory of
Dilute Gases} by Cercignani, Illner, and Pulvirenti
\cite{CIP}(pp. 90-93), who emphasize that the theorem of Lanford
constitutes a {\it validation} of Boltzmann's equation from the
fundamental principles of molecular dynamics.

Let us describe
Lanford's result.

Consider the deterministic dynamics of $n$ hard spheres of radius
$\frac{1}{\sqrt{n}}$.  The phase space is formed by excising
 the points of $(\RR^6)^n$ that represent
configurations in which two or more spheres would overlap.  The set of
all initial configurations that lead to simultaneous collisions of three or
more particles or
 to infinitely many collisions in finite time has measure zero and can
be ignored.  The trajectories through phase space are determined by
the free motion of the spheres between collisions and the rule
(\ref{collisions}) for binary collisions.
  When a trajectory hits a boundary point of the phase space,
i.e., when a collision occurs, the
trajectory continues from the unique boundary point that the rule of
elastic collision associates to it.  This defines the
deterministic dynamics of a dilute gas of $n$ hard spheres of radius
$\frac{1}{\sqrt{n}}$.  Increasing
$n$ increases the number of particles but decreases the
density, whence the term ``dilute limit.''

 Lanford's theorem states (roughly) that
there exists $\tau > 0$ on the order of the mean free time such
that, if the initial $n$-particle densities are $f_0$-chaotic {\it in a very
strong sense }, then the densities at a later time $ t \le \tau $ are
$f_t$ chaotic, where $f_t$ is a mild solution of the Boltzmann equation
with initial
data $f_0$.  The hypotheses on the initial data are that
 the $k$-marginals of the symmetric $n$-particle distributions are
absolutely continuous with continuous densities, and those densities
satisfy a growth bound depending on $k$
and {\it converge uniformly on compact sets to} $f_0^{\otimes k}$ in
the dilute limit $n \longrightarrow \infty$.  This hypothesis on the
initial data is stronger than mere chaos, and Lanford's theorem
asserts that such strong initial chaos is propagated.  The
$n$-particle densities at a later time will be chaotic, says Lanford's
theorem, but typically not chaotic in the same strong sense as were
the initial densities.
This ``loss of convergence quality'' is what permits the Boltzmann
equation to be irreversible even though it is derived from reversible dynamics
\cite[p. 97]{CIP}.

Lanford's theorem says that chaos is
propagated, but only if the initial densities converge uniformly on
compact sets, et cetera.  This kind of propagation of chaos differs
from propagation of chaos as defined in this dissertation; it has to
do with subtler properties of uniform and pointwise convergence of
densities rather than simple weak convergence of distributions.  We remark
that the hard sphere gases of Grad and Lanford do {\it not} propagate chaos
in our sense, nor do they satisfy the conclusions of our theorems
about families of Markov processes that propagate chaos.

\section{Plasmas and Stellar Systems}
\label{Plasmas}

This section contains an account of the propagation of chaos for the
Vlasov equation.

\subsection{Vlasov's equation}
\label{Vlasov's equation}

Vlasov's equation \cite{BH, Spohn} is another important equation of kinetic
theory.  It
governs the density in position-velocity space of
particles that interact (without colliding) through long-range
forces such as the electric forces between ions in a
plasma or the gravitational attraction between stars in a galaxy.

Suppose, for simplicity, that all particles in the system are of the
same species, each having mass $m$, and let $F(\bx)$ denote the force that
a particle at
the origin would exert on a particle at $\bx$.  For example,  the force
$F(\bx)$ is proportional
to $\bx / \| \bx \|^3$ if
the particles are electrons, and proportional to $- m^2 \bx / \| \bx \|^3$
if the particles are stars.  If $f(\bx,\bv,t)d\bx d\bv$ denotes the number
of particles per unit volume near $(\bx,\bv)$ at time $t$, we find
that the net force on a particle at $\bx$ is
\begin{equation}
    F_f(\bx) := \int_{\RR^3} F(\bx - \bxp) f(\bxp,\bvp,t)d\bxp d\bvp  .
\label{meanfield}
\end{equation}
The particle density $f(\bx,\bv,t)$ changes through the motion of
particles subject to the force field $F_f(\bx)$.

Vlasov's equation for the density is
\begin{equation}
\label{FirstVlasov}
   \ppt f(\bx,\bv,t) = -\bv \cdot \nabla_{\bx}f(\bx,\bv,t) -
   \frac{1}{m} F_f(\bx) \cdot \nabla_{\bv} f(\bx,\bv,t) ,
\end{equation}
where the net force field $F_f(\bx)$, defined in equation
(\ref{meanfield}), depends on the particle density $f$ itself.
This is just an advection equation for the flow on
$(\bx,\bv)$-space given by the time dependent flux $(\bv, F(\bx))$,
with the requirement that $F$ equals $F_f$, i.e., the flux function at
time $t$ is determined through (\ref{meanfield}) by the solution at time $t$
of the advection equation itself.

The preceding is a heuristic derivation of Vlasov's equation from the
smoothed dynamics of a large but fixed number of particles.   Vlasov's
equation may be derived rigorously from the true dynamics of
interacting particle systems, in the limit of infinite particle
number.  This rigorous derivation is the content of the theorems,
stated in the next section,
on the propagation of chaos for Vlasov and McKean-Vlasov equations.

Propagation of chaos clarifies the relationship
between the Vlasov equation and the dynamics of gravitational systems
and plasmas.

Imagine $n$ particles of mass $\oon$ following the classical $n$-body
evolution.   As the
number of particles tends to infinity and the initial distribution
of particles approaches a distribution $f(\bx,\bv,0)d\bx d\bv$ of mass,
Vlasov's equation is an increasingly correct
description of the evolution of the mass density.  The mass density
follows equation (\ref{FirstVlasov}) with $m$ set to $1$ and with
$F(\bx)$
redefined as the gravitational
force on a test particle of {\it unit} mass at $\bx$ due to a particle of
{\it unit} mass at the origin.

It is a little tougher
to obtain a macroscopic equation for the density of charge in the
limit of infinitely many ions.
A $k$-fold increase of the number of electrons (say) in a
plasma increases the forces by a factor of $k^2$ and the system
becomes too energetic in the limit $n \longrightarrow \infty$.
Vlasov dynamics can only result from proper scaling of mass and/or time.  One
possibility is to imagine $n$ electrons of mass $n$ each.  In the
limit $n \longrightarrow \infty$, the density of charge in
position-velocity space satisfies equation (\ref{FirstVlasov}),
{\it mutatis mutandis}.  An
alternative scaling is found in \cite{BH}:
Consider the dynamics of
$n$ ions of mass $\oon$ each.  As $n$-tends to infinity, and time is
slowed as $\oon$, one obtains a Vlasov equation for the density in
position-{\it momentum} space.

\subsection{ McKean-Vlasov particle systems}
\label{McKean-Vlasov particle systems}

We begin with a statement of the propagation of chaos for Vlasov's
equation.  It is assumed that the interparticle force is bounded and
globally Lipschitz, an assumption that excludes the physical
inverse-square forces of gravitational systems and plasmas.   One way
around this difficulty is to assume that the system is so dilute that
particles never get too close to one another.   The
 interparticle force could then be replaced with one  without the
singularity at zero distance that is still
inversely proportional to the square of
the distance between particles when that distance is not too
small.

Let $F:\RR^6 \longrightarrow \RR^6$ be bounded and Lipschitz.
   For each $n$, define a
deterministic $n$-particle process in $\RR^6$ by
the following system of ordinary differential equations (ODEs):
\begin{eqnarray}
\label{VlasovSystem}
      \ddt \bx^n_i(t) &=& \bv^n_i(t)   \nonumber \\
      \ddt \bv^n_i(t) &=& \oon\sum_{j=1}^n  F(\bx^n_i - \bx^n_j)
     \nonumber  \\
\end{eqnarray}
for $i = 1,2,\ldots,n$.  Braun and Hepp~\cite{BH} prove that if the initial
conditions
\[
   \bx^n_1(0), \bv^n_1(0),  \bx^n_2(0), \bv^n_2(0), \ldots,  \bx^n_n(0),
\bv^n_n(0)
\]
are such that
\[
     \oon \sum_{i=1}^n \delta_{( \bx^n_i(0), \bv^n_i(0)) }
     \longrightarrow \mu_0 \in \P(\RR^6),
\]
then, for each $t>0$,
\[
        \oon \sum_{i=1}^n \delta_{( \bx^n_i(t), \bv^n_i(t))}
     \longrightarrow \mu_t ,
\]
where $\mu_t \in \P(\RR^6)$ is the weak solution at time $t$ of the
Vlasov equation
\begin{eqnarray}
\label{VlasovEq}
   \ppt f(\bx,\bv,t) &=& -\bv \cdot \nabla_{\bx}f(\bx,\bv,t) -
   F_f(\bx) \cdot \nabla_{\bv} f(\bx,\bv,t)  \nonumber \\
    F_f(\bx) &=& \int_{\RR^3} F(\bx - \bxp) f(\bxp,\bvp,t)d\bxp d\bvp
    \nonumber \\
    \mu_0 &=& f(\bx,\bv,0)d\bx d\bv.  \nonumber \\
\end{eqnarray}

Thanks to our Corollary~\ref{deterministic}, this theorem of Braun and Hepp
implies
that the family of $n$-particle processes (\ref{VlasovSystem}) propagates
chaos.  The fact that the result of Braun and Hepp implies the propagation
of chaos
is also noted in \cite[p. 99]{Pulvirenti}.

The deterministic particle systems (\ref{VlasovSystem}) may be generalized to
interacting {\it diffusions}.  A diffusion is a Markov process with
continuous trajectories, like the solution of a stochastic
differential equation.  McKean~\cite{McK68} initiated the study of
propagation of chaos for diffusions and what is now called the
McKean-Vlasov equation.

Let  ${\bf v}:\RR^d \times \RR^d \longrightarrow \RR^d$ and $\s:\RR^d
\times \RR^d \longrightarrow \RR$
be bounded and globally Lipschitz.
 For each $n$,  consider the system of $n$ stochastic differential
equations (SDEs)
\begin{equation}
   dX_i^n = \left\{ {1 \over n}\sum_{j=1}^n {\bf v}(X^n_i,X^n_j)
   \right\} dt
	    + \left\{ {1 \over n}\sum_{j=1}^n \s (X^n_i,X^n_j) \right\}
dW_i ,
\label{McKeanSystem}
\end{equation}
for random vectors $X^n_1,X^n_2,\ldots,X^n_n$ in ${\mathbb R}^d$.
The Wiener processes
\[  W_1,\quad W_2,\quad W_3,\quad \ldots \]
are taken to be independent of one another and of the random initial conditions
\[
     X_1^n(0), \quad X_2^n(0), \quad \ldots, \quad X_n^n(0)   .
\]
Each system of SDEs has a unique solution and defines a Markov
transition function
\[
     K_n(\bx,d\by,t):(\RR^d)^n \times \B_{(\RR^d)^n} \times
     [0,\infty) \longrightarrow [0,1]
\]
by
\[
    \int_{\RR^d} \phi(\by) K_n(\bx,d\by,t) :=
\EE^{\bx}[\phi(X_1^n(t),\ldots,X_n^n(t))].
\]
In other words, for fixed $t \ge 0$ and $\bx \in (\RR^d)^n $,
$K_n(\bx,\cdot,t)$
 is the distribution of the position at time $t$ of random trajectory
\[
  {\bf X}^n \equiv
   (X^n_1(t), X^n_2(t), \ldots, X^n_n(t)) \in
      \left[ C( [0,\infty), \RR^d) \right]^n
\]
that started at ${\bf X}^n(0) = \bx$.

McKean \cite{McK66,McK68} proves that if the particles are initially
stochastically independent but with a common distribution $\mu_0$,
then the sequence of $n$-particle joint distributions at time $t$ is
$\mu_t$-chaotic, $\mu_t$ being the (weak) solution at time $t$ of the
nonlinear McKean-Vlasov equation
\begin{eqnarray}
\label{McKeanVlasvEq}
      \ppt f_t(\bx) &=& -\nabla\cdot\left[ V_f(\bx) f_t(\bx)  \right] +
      \frac{1}{2} \Delta\left[D_f(\bx) f_t(\bx)\right]   \nonumber  \\
      V_f(\bx) &=& \int_{\RR^d} \bv(\bx,\bxp) f_t(\bxp)d\bxp   \nonumber \\
      D_f(\bx) &=& \left( \int_{\RR^d} \s(\bx,\bxp) f_t(\bxp)d\bxp
\right)^2  \nonumber \\
      f_0(\bx)d\bx &=& \mu_0 . \nonumber  \\
\end{eqnarray}

                 McKean's result includes that of
Braun and Hepp: when $\s \equiv 0$ there is no diffusion and the
system of SDEs (\ref{McKeanSystem}) becomes a Vlasov system of ODEs
like (\ref{VlasovSystem}).  Braun and Hepp seem unaware, in their
paper of 1977, of McKean's important work of 1966.  They use a different
method to prove the propagation of chaos for Vlasov systems.
Though they only treat the deterministic (Vlasov) case,
their method can be generalized to prove that
interacting (McKean-Vlasov) diffusions also propagate chaos.

McKean really proves much more than the propagation of chaos.
Suppose the initial positions $ X_1^n(0), \ldots,  X_n^n(0)$ for the
$n$-particle systems
are taken to be the first $n$ terms of a sequence
$    Z_1, Z_2, Z_3, \ldots $ of independent and $\mu_0$ distributed
random variables.  McKean proves that $X_i^n(t)$, the random position of the
$i^{th}$ particle at time $t$, converges in mean square to
$X_i^{\infty}(t)$ as $n$ tends to infinity.  The $X_i^{\infty}$
are independent and identically distributed.   $X_1^{\infty}(t)$
is sometimes called the {\it nonlinear process} and it satisfies the
SDE
\begin{eqnarray*}
         dX & = & \left\{\int_{\RR^d}\bv(X,y)\mu_t(dy)\right\}dt +
                    \left\{\int_{\RR^d}\sigma(X,y)\mu_t(dy)\right\}dW_1 \\
      \mu_t & = & \mathrm{Law}(X(t))    \\
\end{eqnarray*}
with  $  X(0) = Z_1  $.

\chapter{Chaos and Weak Convergence}
\label{Chaos and Weak Convergence}
Chaos of a sequence of symmetric measures is equivalent to weak convergence
of certain probability measures.  This observation, due to Sznitman
and Tanaka, is the subject of this chapter.
First,  in Section \ref{Background}, the theory of weak convergence
of probability measures is reviewed.  The theorem of Sznitman and
Tanaka is proved in Section \ref{Theorem of Sznitman and Tanaka}.  We
examine this
equivalence in the simplest context of finite probability spaces in Section
\ref{Chaos on Finite Sets}.

This chapter ends with Theorem \ref{SpecEntThm}: on a finite space, chaos
implies convergence of specific entropy to the entropy of the
single-particle distribution.

\section{Background}
\label{Background}

Let $X$ be a set and $\F$ a class of subsets of $X$ that contains
the empty set and is closed under complementation and countable unions.
$(X,\F)$ is called a {\it  measurable space}, and the sets in $\F$
are called {\it measurable}. A {\it probability
measure} or {\it law} on $(X,\F)$ is a countably additive, nonnegative
function
\[
    P:\F \longrightarrow [0,1]
\]
satisfying $P(X) = 1$.
The measure $P(F)$ of a set $F \in \F$ is the {\it probability of} $F$.
Countable additivity
requires the probability of a union of a sequence of disjoint
measurable sets to  equal the sum of their probabilities.  The
simplest probability measure is a {\it point mass} at a point $x \in
X$, denoted $\delta(x)$ or $\delta_x$; $\delta_x(F)$ equals one if $x \in
F$, otherwise it
equals zero.

 Let $(X,\F)$
and $(Y,\G)$ be measurable spaces, and let $h:X \longrightarrow Y$
be {\it measurable}, i.e., $h^{-1}(G) \in \F$ whenever $G \in \G$,
where $h^{-1}(G)$ is the inverse image of $G$ under $h$.
Any probability measure $P$ on $X$ induces a probability measure $P
\circ h^{-1} $ on $Y$
via $h$.  The probability measure {\it induced by} $h$ is defined for $G
\in \G$ by
\[
     \left( P \circ h^{-1}\right)(G) := P( h^{-1}(G))     .
\]
This definition implies that for any integrable function $\phi$ on $(Y,\G,P
\circ h^{-1})$,
\[
     \int_Y \phi(y) P \circ h^{-1}(dy) = \int_X \phi(h(x)) P(dx) .
\]

Now let $(X,\T)$ be a Hausdorff topological space with topology $\T$.  The
{\it Borel $\s$-algebra}, $\B$, is the smallest $\s$-algebra containing $\T$.
The Borel algebra is thus the smallest $\s$-algebra with respect to
which any function continuous on $(X,\T)$
is measurable.    The set of probability
measures on $(X,\B)$ is denoted ${\P}(X)$.

We often call probability measures simply ``laws.''

Let $C_b(X)$ denote the continuous and bounded real-valued functions
on $(X,\T)$.
The set of laws $\P(X)$ is endowed with the weakest topology rendering
continuous the maps
\[
      P  \in \P(X) \longmapsto \int_X g(x) P(dx) \in \RR ,
\]
for all $g \in C_b(X)$.  This is known as the weak topology
on $\P(X)$.   A net of laws $\{ P_{\beta} \}$ in $\P(X)$ converges to $P$
in the weak
topology if and only if the nets $\{ \int g P_{\beta} \}$ converge
to $\int g P$ for all $g \in C_b(X)$.

We consider exclusively the case that $X$ is homeomorphic to a separable
metric space
$(S,d_S)$, so that we may use certain results of the theory of weak
convergence.
The theory of weak convergence of laws is customarily expounded for
laws on separable metric spaces, and
especially {\it complete} and separable metric spaces, because of the
influence of Prohorov's original study \cite{Prohorov} of 1956.
Around the same time, Le Cam \cite{Le Cam} developed the theory of
of weak convergence of laws on completely regular topological
spaces.

For separable metric spaces $(S,d_S)$, the weak topology on $\P(S)$ is
metrizable.  Two metrics on $\P(S)$ that generate the weak topology
are the L\'evy-Prohorov distance $LP$ and Dudley's distance
$BL^*$.   The Dudley distance between two laws $\mu,\nu \in \P(S)$
is
\[
     BL^*(\mu,\nu) := \sup_{g \in BL_1} \left\{ \left| \int_S g(s)
     \mu(ds) - \int_S g(s) \nu(ds)  \right| \right\}  ,
\]
where $g$ ranges over the class $BL_1$ of bounded Lipschitz functions
from $S$ to $\RR$ defined as
\[
      BL_1 := \left\{ g(s): \sup_{s \in
      S}\left\{ |g(s)| \right\} + \sup_{s \ne t \in S}\left\{
      |g(s) - g(t)| / d_S(s,t) \right\}   \le 1  \right\}  .
\]
The L\'evy-Prohorov distance between $\mu$ and $\nu$ is
\[
      LP(\mu,\nu) := \inf \left\{ \delta > 0 : \nu(B) \le
      \mu(B^{+\delta}) +
      \delta  \  \mathrm{ for \  all \ closed \ sets } \  B  \right\}  ,
\]
where $B^{+\delta}$ is the set of all points in $S$ that are within
$\delta$ of some point of $B$:
\[
      B^{+\delta} := \left\{ s \in S : d_S(s,B) < \delta \right\}.
\]
These metrics are discussed in Chapter 11 of the textbook {\it Real Analysis
and Probability}, by Dudley \cite{Dudley}.

The general theory of weak convergence in law on
Polish spaces is due to Prohorov.  (A topological space is {\it Polish} if
it is
homeomorphic to a complete, separable metric space.)  Prohorov's theorem
characterizes
compact sets in $\P(X)$ when $X$ is
Polish, much as the Arzel\`a-Ascoli theorem characterizes compactness
in the space of continuous functions on a compact Hausdorff
space.  The Arzel\`a-Ascoli theorem
states that a family of functions on a compact space is relatively compact
(has compact
closure) in the topology of uniform convergence if and only if the
family is equicontinuous and bounded.  Prohorov's theorem states that
a family of laws on a Polish space is relatively compact if
and only if it is {\it tight}.

Tightness is a simple condition:
\begin{definition}
Let $\Sigma \subset \P(X)$ be a family of laws on a topological space.

The family $\Sigma$ is {\bf tight} if for each $\epsilon > 0$ there
exists a compact $K_{\epsilon} \subset X$ such that
\[
      \sup_{\mu \in \Sigma}\left\{ \mu(X \setminus K_{\epsilon})
      \right\} < \epsilon.
\]
\end{definition}
Tightness implies relative compactness, and the conditions are
equivalent in separable, topologically complete spaces:
\begin{theorem}[Prohorov]
\label{Prohorov}
Suppose $(X,\T)$ is homeomorphic to a separable metric space.  Then, if
$\Sigma
\subset \P(X)$ is tight, its closure is compact in $\P(X)$.

If $(X,\T)$ is Polish (homeomorphic to a complete, separable metric
space) then $\Sigma \subset \P(X)$ is tight if and only if its closure is
compact in $\P(X)$.
\end{theorem}

It follows from Prohorov's theorem that $\P(X)$ is Polish if $X$ is
Polish \cite{Dudley}.  This is important to us since the proof of our
main theorem requires the application of Prohorov's
theorem to $\P(X)$.

{\it Convergence of Probability Measures} by Patrick Billingsley
\cite{Billingsley} is a charming classic monograph on the theory of weak
convergence of laws and its applications.  Unfortunately, this text is
missing some essential material, especially the
metric approach to weak convergence.
It is well complemented by the material in \cite{Dudley}.

\section{The Theorem of Sznitman and Tanaka}
\label{Theorem of Sznitman and Tanaka}

Let $(S,d_S)$ be a separable metric space with Borel algebra $\B_S$.
Let $S^n$ denote the n-fold product of $S$ with itself;
\[
      S^n := \left\{ (s_1, s_2, \ldots, s_n): s_i \in S \quad \mathrm{ for}
       \quad i=1,2,\ldots,n \right\}  .
\]
$S^n$ is itself metrizable in a variety of equivalent ways that all
generate the same topology and the same Borel algebra $\B_{S^n}$.

The {\it marginal} of a law $\rho_n \in (S^n)$ on the first
$k$-coordinates ($k\le n)$ is the law $\rho_n^{(k)} \in \P(S^k)$ induced by
the projection
\[
     (s_1, s_2, \ldots, s_n) \longmapsto  (s_1,s_2, \ldots, s_k).
\]
Equivalently,
\[
     \rho_n^{(k)}(B_1,B_2, \ldots, B_k) = \rho_n(B_1,B_2, \ldots,
     B_k,S,S,\ldots,S)  ,
\]
for all $B_1,B_2, \ldots, B_k \in \B_S$.
If $\rho \in \P(S)$, the {\it product law} $\rho^{\otimes n} \in \P(S^n)$
is the law
$\rho(ds_1)\rho(ds_2)\cdots \rho(ds_n)$.  Note that $\left( \rho^{\otimes
n} \right)^{(k)} = \rho^{\otimes k}$.

Let $\Pi_n$ denote the set of permutations of $\{1,2,\ldots,n\}$.
The permutations $\Pi_n$ act on $S^n$ by permuting coordinates:
the map $\pi\cdot:S^n \longrightarrow S^n$ is
\[
    \pi \cdot (s_1, s_2, \ldots, s_n)  := (s_{\pi(1)},s_{\pi(2)},
    \ldots,s_{\pi(n)})  .
\]
If $E$ is any subset of $S^n$, define
\[
     \pi\cdot E = \{ \pi\cdot\bs: \bs \in E \}  .
\]
A law $\rho$ on $S^n$ is {\it symmetric} if $\rho(\pi\cdot B) = \rho(B)$
for all $\pi \in \Pi_n$ and all $B \in \B_{S^n}$.  Products $\rho^{\otimes
n}$ are symmetric,
for example.  The {\it
symmetrization} ${\widetilde \rho}$ of a law $\rho \in \P(S^n)$ is the
symmetric law such that
\[
    {\widetilde \rho}(B) := \frac{1}{n!}\sum_{\pi \in \Pi_n} \rho(\pi\cdot
    B),
\]
for all $B \in \B_{S^n}$.

\begin{definition}[Kac, 1954]
\label{chaos}
Let $(S,d_S)$ be a separable metric space.  Let $\rho$ be a law on $S$, and
for $n = 1,2 ,\ldots$,
let $\rho_n$ be a symmetric law on $S^n$.

The sequence
$\{ \rho_n \}$ is $\rho$-{\bf chaotic} if, for each natural number $k$
and each choice
\[
     \phi_1(s),\quad \phi_2(s),\quad \ldots, \quad\phi_k(s)
\]
 of $k$ bounded and continuous functions on $S$,
\begin{equation}
\label{KacCondition}
\lim_{n \rightarrow \infty} \int_{S^n} \phi_1(s_1) \phi_2(s_2) \cdots
\phi_k(s_k) \rho_n(ds_1 ds_2\ldots ds_n)    =
\prod_{i=1}^k \int_S \phi_i(s) \rho(ds).
\end{equation}
\end{definition}

In case $S$ is Polish, condition \ref{KacCondition} implies
the weak convergence of the marginals to products $\rho^{\otimes k}$,
because the class of functions of the form
\begin{equation}
\label{products}
     \phi_1(x_1) \phi_2(x_2) \cdots \phi_k(x_k);\qquad \phi_1,
     \ldots,\phi_k \in C_b(S)
\end{equation}
is a {\it convergence determining class} for $\P(S^k)$\cite{EthKur}.
Condition (\ref{KacCondition}) shows that the sequence of the marginals
$\rho_n^{(k)}$ converges to $\rho^{\otimes k}$ weakly {\it for
functions of the form} \ref{products} , hence it converges weakly.
 Thus, {\it if $S$ is Polish}, a sequence $\{ \rho_n \}$ of symmetric
laws on $S^n$ is $\rho$-chaotic if and only if
\[
    \lim_{n \rightarrow \infty} \rho_n^{(k)} = \rho^{\otimes k}  ,
\]
for any natural number $k$.

It turns out, however, that $S$ does not need to be Polish.
It will be seen from the proof of the next theorem that condition
(\ref{KacCondition})
implies the convergence of the marginals $\rho_n^{(k)}$ {\it even if
$S$ is not Polish}, but only separable.

The following theorem of Sznitman and Tanaka
states that a sequence of symmetric laws is chaotic if and only if
the induced sequence of laws of the random empirical measures converges to
a point mass.  Let
\begin{equation}
\label{e's}
        \e_n((s_1,s_2,\ldots,s_n)) := \oon\sum_{i=1}^n \delta(s_i)
\end{equation}
define a map from $S^n$ to $\P(S)$.  These maps are measurable for
each $n$, and $\e_n(\pi \cdot\bs)  = \e_n(\bs)$ for all $\bs \in S^n, \pi
\in \Pi_n$.

\begin{theorem}[Sznitman, Tanaka]
\label{SznitmanTanaka}
$\{ \rho_n \}$ is $\rho$-chaotic if and only if
\begin{equation}
\label{SznitmanCondition}
    \rho_n \circ \e_n^{-1} \longrightarrow \delta(\rho)
\end{equation}
in $\P(\P(S))$.
\end{theorem}

\noindent{\bf Proof}:

Suppose $\{ \rho_n \}$ is $\rho$-chaotic.

A sequence of laws $\{ \mu_n \}$ on a completely
regular topological space $X$ converges to
$\delta(x) \in \P(X)$ if and only if for each neighborhood $N$ of $x$
\begin{equation}
\label{xxx1}
     \lim_{n \longrightarrow \infty} \mu_n(X \setminus N) = 0
\end{equation}
for each neighborhood $N$ of $x$.
Therefore, to prove the convergence of $\rho_n \circ \e_n^{-1}$ to
$\delta(\rho)$ in $\P(\P(S))$ it suffices to verify (\ref{xxx1}) on a
subbase of neighborhoods of $ \rho \in \P(S)$.  The class of sets of the form
\[
      \left\{ \nu \in \P(S) : \left| \int_S g(s) \nu(ds) -  \int_S
      g(s) \rho(ds) \right| < \epsilon \right\}; \quad \epsilon >
      0, g \in C_b(S)
\]
is a neighborhood subbase at $\delta(\rho)$, so it suffices to show
that
\begin{equation}
\label{xxx2}
      \rho_n \circ \e_n^{-1} \left( \left \{ \nu  : \left| \int_S g(s)
      \nu(ds) -  \int_S g(s) \rho(ds) \right| \ge \epsilon \right \}
      \right)  \longrightarrow 0.
\end{equation}

Writing $\int_S g(s) \nu(ds)$ as $<g,\nu>$, we calculate
\begin{eqnarray*}
      && \int_{S^n}  \left| <g ,\e_n(\bs)> - <g,\rho> \right|^2
       \rho_n(d\bs)   \\
       & = & \int_{S^n} \left( \oon\sum_{i=1}^n g(s_i) - <g,\rho>
       \right)^2
       \rho_n(ds_1ds_2\cdots ds_n)   \\
       & = &  \frac{1}{n^2} \sum_{i,j=1}^n \int_{S^n}\left( g(s_i) -
       <g,\rho> \right)\left( g(s_j) - <g,\rho> \right) \rho_n(d\bs)   \\
       & = & \oon\int_S \left( g(s) - <g,\rho> \right)^2
       \rho_n^{(1)}(ds)  \\
       &\ & +  \frac{n-1}{n} \int_{S\times S}\left( g(s_1) -
       <g,\rho> \right)\left( g(s_2) - <g,\rho> \right)
       \rho_n^{(2)}(ds_1ds_2)  ,   \\
\end{eqnarray*}
the last equality by the symmetry of $\rho_n$.  Thus condition
(\ref{KacCondition}) for $k = 1,2$ implies that
\[
   \int_{S^n}  \left| <g ,\e_n(\bs)> - <g,\rho> \right|^2
       \rho_n(d\bs)   \longrightarrow 0 ,
\]
and hence that (\ref{xxx2}) holds.  Condition
(\ref{KacCondition}) thus implies condition (\ref{SznitmanCondition}).

Now suppose that $\rho_n \circ \e_n^{-1} $  tends to $
\delta(\rho)$.

For natural numbers $k \le n$, let $\J_{n:k}$ and $\I_{n:k}$ denote
respectively
the set of all maps and the set of
injections from $\{1,2,\ldots,k\}$ into $\{1,2,\ldots,n\}$.  Define
the map $\e_{n:k}$ from $S^n$ to $\P(S^k)$ by
\begin{equation}
\label{e:k's}
     \e_{n:k}((s_1,s_2,\ldots,s_n)) :=
     \frac{(n-k)!}{n!}\sum_{i \in
     \I_{n:k}}\delta_{(s_{i(1)},\ldots,s_{i(k)})}  .
\end{equation}
$\e_{n:k}(\bs) $ is the empirical measure of $k$-tuples of
coordinates of $\bs$, sampled without replacement.  Define also
\begin{equation}
\label{varthetas}
        \vartheta_{n:k}(\bs) := \e_n(\bs)^{\otimes k} =
\frac{1}{n^k}\sum_{j \in
               \J_{n:k}}\delta_{(s_{j(1)},\ldots,s_{j(k)})} \ ,
\end{equation}
the empirical measure of {\it all} $k$-tuples from $\bs$.
 When $n >> k$, these two empirical measures are close in total variation
 (TV) and {\it a fortiori} in Dudley's distance on $\P(S)$:
\begin{eqnarray*}
       && BL^*\left( \e_{n:k}(\bs) , \e_n(\bs)^{\otimes k}
               \right)  \\
       &\le& \left\| \frac{(n-k)!}{n!}\sum_{i \in
       \I_{n:k}}\delta_{(s_{i(1)},\ldots,s_{i(k)})}
               - \frac{1}{n^k}\sum_{j \in
               \J_{n:k}}\delta_{(s_{j(1)},\ldots,s_{j(k)})}
               \right\|_{TV} \\
       &\le& 2\left(1 - \frac{n!}{n^k(n-k)!}\right) \ .  \\
\end{eqnarray*}
Since this bound is uniform in $\bs$, it follows that  $ \rho_n\circ
\e_{n:k}^{-1}$ is near $ \rho_n \circ
\vartheta_{n:k}^{-1}$ in $\P(\P(S^k))$.  In fact,  both
 the L\'evy-Prohorov  and the Dudley distances between the two laws
are bounded above:
\begin{eqnarray}
\label{distances}
       BL^*\left(\rho_n\circ \e_{n:k}^{-1}, \rho_n \circ \vartheta_{n:k}^{-1}
       \right )  & \le & 2\left(1 - \frac{n!}{n^k(n-k)!}\right)  \nonumber
\\
      \mathrm{and} \qquad \qquad \qquad \quad
       LP\left(\rho_n\circ \e_{n:k}^{-1}, \rho_n \circ \vartheta_{n:k}^{-1}
       \right )  & \le &  2\left(1 - \frac{n!}{n^k(n-k)!}\right) \quad.
\nonumber \\
\end{eqnarray}

Condition (\ref{SznitmanCondition}) and definition (\ref{varthetas})
imply that $\rho_n \circ \vartheta_{n:k}^{-1}$ converges to
$\delta(\rho^{\otimes k})$ in $\P(\P(S^k))$.   By (\ref{distances}),
$\rho_n\circ \e_{n:k}^{-1}$ converges to $\delta(\rho^{\otimes k})$ as well.

Now, if $\phi \in C_b(S^k)$,
\begin{eqnarray*}
     &&  \lim_{n \rightarrow \infty}\int_{S^n}
     \phi(s_1,s_2,\ldots,s_k)\rho_n(d\bs)     \\
     & = &  \lim_{n \rightarrow \infty}   \int_{S^n} \left\{
\frac{(n-k)!}{n!}\sum_{i \in
     \I_{n:k}} \phi(s_{i(1)},\ldots,s_{i(k)}) \right\} \rho_n(d\bs)  \\
     & = &   \lim_{n \rightarrow \infty}  \int_{S^n}   <\phi,
\e_{n:k}(\bs)>  \rho_n(d\bs) \\
     & = & \lim_{n \rightarrow \infty} \int_{\P(S^k)} <\phi,\mu>
\rho_n\circ \e_{n:k}^{-1}(d\mu) \\
     & = & \lim_{n \rightarrow \infty} \int_{\P(S^k)} <\phi,\mu> \rho_n\circ
     \vartheta_{n:k}^{-1}(d\mu) \\
     & = & \int_{\P(S^k) }<\phi,\mu> \delta(\rho^{\otimes k})(d\mu)    \\
     & = & \int_{S^k} \phi(s_1,s_2,\ldots,s_k) \rho(ds_1)\cdots\rho(ds_k).  \\
\end{eqnarray*}
Thus, condition (\ref{SznitmanCondition}) implies (\ref{KacCondition}).
$\blacksquare$

The preceding arguments have actually proved the following stronger
version of Theorem~\ref{SznitmanTanaka}.

\begin{theorem}
Let $S$ be a separable metric space and for each $n$ let $\rho_n$ be
a symmetric law on $S^n$.

The following are equivalent:
\begin{trivlist}
\item{{\bf  Kac's condition for $k = 2$}:} For all $\phi_1,\phi_2 \in C_b(S)$,
 \begin{equation}
  \lim_{n \rightarrow \infty} \int_{S^n} \phi_1(s_1) \phi_2(s_2)
  \rho_n(d\bs)  = \int_S \phi_1(s) \rho(ds) \int_S \phi_2(s) \rho(ds) ;
 \end{equation}
\item{{\bf Condition of Sznitman and Tanaka}:} For all natural numbers $k$,
the laws $\rho_n\circ
\e_{n:k}^{-1}$ converge to $\delta(\rho^{\otimes
  k})$ in $\P(\P(S))$ as $n$ tends to infinity,
  where $\e_{n:k}$ is the empirical measure defined in (\ref{e:k's}) ;
\item{{\bf Weak convergence of marginals}:} For all $k$,
  the marginals $\rho_n^{(k)}$ converge weakly to $\rho^{\otimes k}$
  as $n$ tends to infinity.
\end{trivlist}

\end{theorem}

\section{Chaos on Finite Sets}
\label{Chaos on Finite Sets}

Throughout this section, let $S = \{ s_1,s_2,\ldots,s_k \}$ be a finite set.

For each natural number $n$, let
\[ \rho_n(x_1,x_2,\ldots,x_n) \]
be a symmetric law on
$ S^n = S\times S \times \cdots \times S $.
Because of its symmetry, $\rho_n$
is entirely determined by the probability function
\begin{equation}
   P_n(j_1,j_2,\ldots,j_k) ; \quad \sum_{i=1}^k j_i = n
\end{equation}
that gives the probability there are $j_1$ coordinates equal to $s_1$,
$j_2$ coordinates equal to $s_2$, and so on.
The probability of  $(x_1,x_2,\ldots,x_n)\in S^n$ is
\begin{equation}
    \rho_n(x_1,x_2,\ldots,x_n) =
      {P_n(j_1,j_2,\ldots,j_k) / {n! \over j_1! \cdots j_k!}}
      \quad ,
\label{rho'n'P}
\end{equation}
where $j_i(x_1,x_2,\ldots,x_n)$ is the number of coordinates of
$(x_1,x_2,\ldots,x_n)$
that equal $s_i$.

Let $\Delta_{k-1}$ denote the unit simplex in ${\mathbb R}^k$:
\[
\Delta_{k-1} = \left\{ (q_1,q_2,\ldots,q_k):
			    \sum_{i=1}^k q_i = 1 ,  q_i \ge 0
			    \right\}.
\]
Given $\rho_n$, define a law
$\mu_n$ on $\Delta_{k-1}$ by
\begin{equation}
\label{mu}
      \mu_n := \sum_{{\bf j} }
      P_n({\bf j})
      \delta (  {\bf j}/n) \quad,
\end{equation}
where ${\bf j}$ ranges over $k$-tuples of nonnegative integers that
sum to $n$.

Finally, let ${\bf p} = (p_1,p_2,\ldots,p_k)$ be a point of
$\Delta_{k-1}$, and let $ p$ denote the law on $S$ given by $p(s_i) =
p_i$.  With these definitions and notations, we can formulate simpler
versions of Definition \ref{chaos} and Theorem \ref{SznitmanTanaka}
for finite probability spaces:

\begin{definition}[Chaos for Finite State Spaces]

The sequence $\left\{ \rho_n \right\}$ is {\bf $p$-chaotic} if for each
natural number  $m$
and each $(z_1,z_2,\ldots,z_m)$ in $S^m$,
\[
  \lim_{n \rightarrow \infty}
  \sum_{x_1,\ldots,x_{n-m} \in S}
  \rho_n(z_1,z_2,\ldots,z_m,x_1,x_2,\ldots,x_{n-m})
  =    \prod_{i=1}^m p(z_i)   .
\]
\end{definition}

\begin{theorem}
\label{FiniteSznitmanTanaka}

The sequence $\{\rho_n \}$ is $ p$-chaotic if and
only if $\mu_n$ converges weakly to $\delta(\bf p)$.

Equivalently, $\{ \rho_n \}$ is $ p$-chaotic if and only if
\begin{equation}
   \lim_{n \rightarrow \infty}\sum_{(j_1,\ldots,j_k)} P_n(j_1,\ldots,j_k)
	 F\left( {j_1 \over n},{j_2 \over n}, \ldots,
	 {j_k \over n}\right) =    F({\bf p}),
\label{SimpleChaos}
\end{equation}
for {\it every} continuous function $F$ on the simplex $\Delta_{k-1}$.

\end{theorem}
\noindent {\bf Proof}:

This is a special case of Theorem \ref{SznitmanTanaka}. $\blacksquare$

\smallskip

Formula (\ref{SimpleChaos}) will be used to prove that, on finite
probability spaces, chaos implies convergence of specific entropy.  We
are borrowing the expression ``specific entropy'' from statistical
mechanics, where it refers to entropy per particle.

For laws $\pi$ and $\mu$ on a measurable space
$(X,{\mathcal F})$, the entropy of $\mu$
relative to $\pi$ is defined to be
\[
      H_{\pi}(\mu) := -\int_X \left[ {d\mu \over d\pi} \right]
      \log \left[ {d\mu \over d\pi} \right] d\pi
\]
if $\mu$ is absolutely continuous relative to $\pi$ with density
$\left[ {d\mu \over d\pi} \right]$, and
to equal $- \infty$ otherwise.

Relative entropy is nonpositive, but might equal $-\infty$.  $H_{\pi}(\mu)$
achieves its maximum of $0$ only when $\mu = \pi$.  If $X$ is
a Polish space,  $H_{\pi}(\mu)$ is a upper
semicontinuous function of $\mu$ relative to the weak topology on
$\P(X)$.  The entropy of a joint law is less
than or equal to the sum of the entropies of its marginals, with
equality only if the joint law is a product measure.  That
is, if $\mu \in \P(X\times X)$ with marginals $\mu_1, \mu_2 \in
\P(X)$, then
\begin{equation}
\label{subadditivity}
     H_{\pi \otimes \pi}(\mu) \le H_{\pi}(\mu_1) + H_{\pi}(\mu_2) ,
\end{equation}
for any {\it reference law} $\pi \in \P(X)$.  The reader is
referred to
\cite[pp. 32-40]{DupuisEllis} for properties of the relative entropy.

Now,
if $\pi$ is a reference law and $\{\rho_n\}$ is a
$p$-chaotic sequence of laws on a general separable metric space
(where chaos has been defined), the subadditivity (\ref{subadditivity}) of
entropy guarantees that
\[
     \limsup_{n \rightarrow \infty} \oon H_{\pi^{\otimes n}}(\rho_n)
     \le   H_{\pi}(p).
\]
The left hand side of this inequality is what we are calling the
specific entropy.
In case $\{\rho_n\}$ is {\it purely} chaotic, i.e., in case $\rho_n =
p^{\otimes n}$ for all $n$, the specific entropy always equals the
entropy of $p$.
At the other extreme, when the symmetric laws of
a $p$-chaotic sequence $\{\rho_n\}$ are not absolutely continuous
relative to the laws $p^{\otimes n}$, the above inequality is strict,
for then
\[
      \lim_{n \rightarrow \infty} \oon H_{p^{\otimes n}}(\rho_n)
      = - \infty <  0 = H_p(p)  .
\]
However, {\it if the space is finite}, one can prove that the specific
entropy of a chaotic sequence does converge:

\begin{theorem}[Specific Entropy Converges]
\label{SpecEntThm}
Let $S = \{s_1,s_2,\ldots,s_k\}$ be a finite set, and for each $n$ let
$\rho_n \in \P(S^n)$ be a
symmetric law.

If the sequence $\{\rho_n\}$ is $ p$-chaotic,
then
\[
  \lim_{n \rightarrow \infty} \left( - {1 \over n}\sum_{{\bf x}\in S^n}
      \rho_n({\bf x}) \log \rho_n({\bf x}) \right) =
      -\sum_{i=1}^k p_i \log p_i.
\]
\end{theorem}

\noindent {\bf Proof:}

By the relationship (\ref{rho'n'P}) between $\rho_n$ and $P_n$,
\[
   - {1 \over n}\sum_{{\bf x}\in S^n}
      \rho_n({\bf x}) \log \rho_n({\bf x})
   =
    - {1 \over n}\sum_{(j_1,\ldots,j_k)} P_n(j_1,\ldots,j_k) \log
    \frac{P_n(j_1,\ldots,j_k)}{ \frac{n!}{j_1!\cdots j_k!}}.
\]
This equals
\begin{equation}
\label{SpecEnt1}
    -\oon \sum_{\bf j} P_n({\bf j}) \log P_n({\bf j})
    + \oon \sum_{\bf j}
	 P_n({\bf j})\log\left(\frac{n!}{j_1!\cdots j_k!}\right),
\end{equation}
abbreviating $(j_1, j_2, \ldots, j_k)$ by ${\bf j}$.
The first addend in (\ref{SpecEnt1}) is $O\left({\log n \over n}\right)$,
since it equals an $n^{th}$ part of the entropy of a probability
function $P_n({\bf j})$ on fewer than $n^k$ points, which entropy
cannot exceed $\log (n^k) = k \log n$.

Using Stirling's approximation
\[
\log j! = j\log j - j +  \epsilon_j
\]
where $0 < \epsilon_j  < 1 + \log j $ ,
and the fact that $n=\sum_{i=1}^k j_i$,
one finds that
\begin{eqnarray*}
    \log\left({n! \over j_1!\cdots j_k!}\right)
    &=&
     n\log n - n + \epsilon_n - \sum_{i=1}^k(j_i\log j_i - j_i)
	 - \sum_{i=1}^k \epsilon_{j_i}                         \\
    &=&
    n\log n - \sum_{i=1}^kj_i\log j_i + O\left(\log n \right)  \\
    &=&
    - \sum_{i=1}^kj_i \log{j_i \over n} + O\left(\log n \right).\\
\end{eqnarray*}
Substituting this into the second term of (\ref{SpecEnt1}) shows
that
\begin{equation}
\label{SpecEnt2}
   - {1 \over n}\sum_{{\bf x}\in S^n}
        \rho_n({\bf x}) \log \rho_n({\bf x})
   =
    - \sum_{\bf j} P_n({\bf j}) \sum_{i=1}^k
     \left({j_i \over n}\right) \log \left({j_i \over n}\right)
     + O\left({\log n \over n}\right) .
\end{equation}

Since $\rho_n$ is $\rho$-chaotic, formula (\ref{SimpleChaos}) of Theorem
\ref{FiniteSznitmanTanaka}
tells us that
\begin{equation}
\label{SpecEnt3}
   \lim_{n \rightarrow \infty} \left\{ - \sum_{\bf j} P_n({\bf j})
\sum_{i=1}^k
     \left({j_i \over n}\right) \log \left({j_i \over n}\right)
     \right\} =   - \sum_{i=1}^k p_i \log p_i.
\end{equation}
By (\ref{SpecEnt2})  and (\ref{SpecEnt3})
the specific entropy converges to the entropy of $p$:
\[
   \lim_{n \rightarrow \infty}
   \left( - {1 \over n}\sum_{{\bf x}\in S^n} \rho_n({\bf x})
   \log \rho_n({\bf x}) \right) = - \sum_{i=1}^k
   p_i \log p_i.    \qquad \blacksquare
\]

\chapter{Propagation of Chaos}
\label{Propagation of Chaos}

This brief chapter is devoted to the proof of the main theorem stated
in Section \ref{Definition of Propagation of Chaos and Statement of
Main Result}.  Definitions are given and
the approach is outlined in Section \ref{Preliminaries}.  Lemmas are
proved in Section \ref{Lemmas} that expedite the proofs of the theorems
of Section \ref{Theorems}.
\section{Preliminaries}
\label{Preliminaries}

Let $(X,{\mathcal F})$ and $(Y,{\mathcal G})$ be two measurable spaces.
A Markov transition function $K(x,E)$ on $X \times {\mathcal G}$ is a
function that satisfies the following two conditions:

\noindent (1) \qquad $K(x,\cdot)$ is a probability
measure on $(Y,{\mathcal G})$ for each $x \in X$, and

\noindent (2) \qquad $K(\cdot, E)$ is a measurable function on
$(X,{\mathcal F})$ for each $E \in {\mathcal G}$.

Whenever $X$ and $Y$ are measurable
spaces and there is no confusion about what their $\sigma$-algebras
are supposed to be, we usually speak of Markov transitions from
$X$ to $Y$ rather than transition {\it functions}.  In particular, if
$S$ and $T$ are metric spaces, a Markov transition from $S$ to $T$ is
a transition function on $S\times \B_T$.

A Markov {\it process} on a state space $(X,{\mathcal F})$ determines a
family, indexed by time, of Markov transitions from $X$ to itself: $\{
K(x,E,t) \}_{t \ge 0}$.
The transitions
satisfy --- in addition to (1) and (2) above --- the Chapman-Kolmogorov
equations
\[
      K(x,E,s+t) = \int_X K(x,dy,s) K(y,E,t); \quad s,t \ge 0,x \in X, E
\in \mathcal{
      F}.
\]

Let $(S,d_S)$ and $(T,d_T)$ be separable metric spaces.  For each $n$, let
$K_n$ be a Markov
transition from $S^n$ to $T^n$.
We assume that the Markov transition function $K_n$
is symmetric in the sense that,
if $\pi$ is a permutation in $\Pi_n$  and $A$ is a Borel
subset of $T^n$,
\begin{equation}
    K_n(\pi\cdot \bs,\pi\cdot A) =  K_n( \bs, A)   .
\label{permute}
\end{equation}

\begin{definition}[Propagation of Chaos]
\label{PoCDef1}

Let $\{K_n\}$ be as above.

The sequence $\{K_n\}$ {\bf propagates chaos} if,
 whenever $\{ \rho_n \}$ is a $\rho$-chaotic sequence of
measures on $S^n$, the measures
\[
  \int_{S_n} K_n({\bf s},\cdot) \rho_n(d{\bf s})
\]
on $T^n$ are $\tau$-chaotic for some $\tau \in \P(T)$.
\end{definition}

 When we say that a family of $n$-particle Markov {\it processes} on a
state space $S$
{\bf  propagates chaos} we mean that, for each fixed time $t > 0$, the
family of associated
 $n$-particle transition functions $\{K_n(\bs,E,t) \}$
 propagates chaos.

Most Markov processes of interest are characterized by their laws on nice
 {\it path spaces}, such as $C([0,\infty), S)$ the space of continuous
paths in $S$, or the space $D([0,\infty), S)$ of right continuous paths in $S$
 having left limits.  For such processes,
 the function that maps a state $s \in S$ to the law of the process
 started at $s$ defines a Markov transition from $S$ to the entire
 path space.
Now, if a sequence of transitions $K_n(s,\cdot)$ from $S^n$
 to the path spaces $C([0,\infty),S^n)$ or $D([0,\infty),S^n)$ propagates
 chaos, then, {\it a fortiori}, it propagates chaotic sequences of
initial laws to chaotic sequences of laws on $S^n$ at any (fixed) later
time.  We have defined the propagation of chaos
for sequences of Markov transitions from $S^n$ to a (possibly) different
space $T^n$, instead of
simply from $S^n$ to itself, with the case where $T$ is path space
especially in mind.
This way, our ensuing study will pertain even
to those families of processes that propagate the chaos of initial
laws to the chaos of laws on the whole path space.

 We are going to prove that a sequence of Markov transitions $\{ K_n
\}_{n=1}^\infty$
propagates chaos if and only if
\[
     \left\{ \Kt_n(\bs_n,\cdot) \right\}_{n=1}^{\infty}
\]
is chaotic whenever $ \bs_n \in S^n $ satisfy $ \e_n(\bs_n) \longrightarrow
p $ in $\P(S)$.
We employ the weak convergence characterization of chaos of
Sznitman and Tanaka and we assume that $S$ is Polish.
To study the {\it propagation} of chaos, we project the transitions $K_n$
from $S^n$ to $T^n$ onto transitions from $\e_n(S^n)$ to $\P(T)$, and
then apply Theorem \ref{SznitmanTanaka},
which projects chaotic sequences of symmetric laws on the
spaces $S^n$ onto convergent sequences of laws on $\P(S)$.

>From now on, the notation $\e_n$
is used both for the map from $S^n$ to $n$-point empirical measures on
$S$ and for the same kind of map on $T^n$.

 Markov transitions $K_n$ from $S^n$ to $T^n$ induce
Markov transition functions $H_n$ from $\e_n(S^n)$ to $\e_n(T^n)$.
The induced transition function can be defined in terms of a Markov
transition $J_n$ from
$\e_n(S^n)$ to $S^n$ which acts as a kind of inverse of $\e_n$.  For fixed
$\zeta \in \e_n(S^n)$, let $J_n(\zeta,\cdot)$ denote the
atomic probability measure on $S^n$ that allots equal probability to
each of the points in $\e^{-1}(\{ \zeta \})$, a set containing at
most $n!$ points.  Putting it another way, $J_n(\zeta,E)$ equals the
proportion of
 points ${\bf s} \in S^n$ such that $\e_n({\bf s}) =
\zeta$ that lie in $E \subset S^n$.
A Markov transition $K_n$ from $S^n$ to $T^n$ induces
a Markov transition $H_n$ from $\e_n(S^n)$ to $\e_n(T^n)$ defined by
\begin{equation}
  H_n(\zeta,G) := \int_{{\bf s} \in S^n} J_n(\zeta, d{\bf s})
        K_n\left( {\bf s}, \e_n^{-1} ( G ) \right)
\label{InducedTransition}
\end{equation}
for $\zeta \in \e_n(S^n)$ and $G$ a measurable subset of  $\e_n(T^n)$.
Note that if $\bs \in S^n$,
\begin{equation}
\label{NeedtoNote}
        H_n(\e_n(\bs), \cdot)  = K_n(\bs,\cdot) \circ \e_n^{-1},
\end{equation}
where the maps $\e_n$ written on the left and right hand sides are,
respectively, the maps from $S^n$ and $T^n$ to empirical measures
in $\P(S)$ and $\P(T)$.

Theorem \ref{SznitmanTanaka} shows that propagation of chaos by a
sequence $K_n$ is equivalent to the following condition on the
induced transitions $H_n$.
\begin{proposition}
\label{PoCDef2}
The sequence of Markov transitions $\{ K_n \}_{n=1}^\infty$ propagates
chaos if
and only if,
whenever $\{ \mu_n \in \P(\e_n(S^n)) \}_{n=1}^\infty$ converges in
$\P(\P(S))$ to
$\delta(p)$, the sequence
\begin{equation}
\label{H sequence}
\left\{ \int_{\e_n(S^n)} H_n(\zeta,\cdot) \mu_n(d\zeta) \right\}_{n=1}^\infty
\end{equation}
converges in $\P(\P(T))$ to $\delta(q)$, for some $q \in \P(T)$.
\end{proposition}
\noindent{\bf Proof}:

$\{ \mu_n \in \P(\e_n(S^n)) \}$ converges to $\delta(p)$
if and only if
\[
   \left\{ \int_{\e_n(S^n)} J_n(\zeta,\cdot) \mu_n(d\zeta)
   \right\}_{n=1}^{\infty}
\]
is chaotic.  Therefore, the sequence of transitions $\{
K_n \}$ propagates chaos if and only if
\begin{equation}
\label{K sequence}
   \left\{  \int_{S^n}
   K_n(\bs,\cdot) \int_{\e_n(S^n)} J_n(\zeta,d\bs) \mu_n(d\zeta)
   \right\}_{n=1}^{\infty}
\end{equation}
is chaotic whenever $ \mu_n \longrightarrow \delta(p)$.

Now, using definition (\ref{InducedTransition}) of the transitions
$H_n$, we find
\begin{eqnarray*}
   &   &
   \left( \int_{S^n}
     K_n(\bs,\cdot) \int_{\e_n(S^n)} J_n(\zeta,d\bs) \mu_n(d\zeta)
     \right) \circ \e_n^{-1}  \\
   & = &
    \int_{\e_n(S^n)} \left( \int_{\bs \in S^n}   J_n(\zeta,d\bs)
K_n(\bs,\cdot)\circ \e_n^{-1} \right)
   \mu_n(d\zeta)  \\
   & = &
    \int_{\e_n(S^n)}  H_n(\zeta,\cdot)
   \mu_n(d\zeta)    .  \\
\end{eqnarray*}

Thus, by Theorem \ref{SznitmanTanaka}, the sequence (\ref{K sequence}) is
chaotic if and only if the sequence
(\ref{H sequence}) converges to a point mass in $\P(\P(T))$.
Therefore, $\{ K_n \}$ propagates chaos if and only if  (\ref{H
sequence}) converges to $\delta(q)$, for some $q \in \P(T)$.
\qquad $\blacksquare$

Proposition \ref{PoCDef2} implies that if a
sequence of Markov transitions $\{ K_n \}$ propagates chaos, then
the sequence $\left\{   H_n(\zeta_n,\cdot)  \right\}$
converges to a point mass whenever  $\{ \zeta_n \in \e_n(S^n) \}$ converges
in $\P(S)$.
That this condition implies propagation of chaos (and is
not just a necessary condition) is equivalent to our main theorem.
 To prove the sufficiency of the condition, we use
Lemma \ref{lemma3} of the next section.

The lemmas of Section \ref{Lemmas} are presented in a general context.  In
Section \ref{Theorems} we apply these lemmas to propagation of
chaos.   In that context
the induced transitions $H_n$, thought of as functions from
$\e_n(S^n)$ to $\P(\P(T))$, behave like the maps $f_n$ of the lemmas.

\section{Lemmas}
\label{Lemmas}

Let $(X,d_X)$ be a metric space, and $D_1 \subset D_2 \subset \cdots $ an
increasing chain of Borel subsets of $X$ whose union is dense in $X$.  For
each
natural number $n$, let $f_n$ be a measurable real-valued function on $D_n$.

Consider the following four conditions on the sequence $\{ f_n
\}_{n=1}^\infty$.  They are listed in order of decreasing
strength.

\begin{trivlist}

\item {[A]} \qquad Whenever $\{ \mu_n \}$ is a weakly convergent sequence of
probability measures on $X$ with $\mu_n$ supported on $D_n$, then
the sequence
\[
\left\{ \int_X f_n(x) \mu_n(dx) \right\}_{n=1}^{\infty}
\]
of real numbers converges as well.

\item {[B]} \qquad Whenever $\{ \mu_n \}$ is a sequence of
probability measures on $X$ that converges weakly to $\delta(x)$ for some
$x \in X$, and $\mu_n$ is supported on $D_n$, then the
sequence $\left\{ \int_X f_n(x) \mu_n(dx) \right\}$ also converges.

\item {[C]} \qquad Whenever $\{ d_n \}$ is a convergent
sequence of points in $X$, with $d_n \in D_n$, then $\{ f_n(d_n) \}$ also
converges.

\item {[D]} \qquad For any compact $K \subset X$, and for any $\epsilon > 0$,
there exists a natural number $N$ such that, whenever $ m \ge n \ge N$ and
$d \in D_n \cap K$, then
\[
 |f_m(d) - f_n(d)| < \epsilon  .
\]
\end{trivlist}

\begin{lemma}
\label{lemma1}
\qquad [A] $\Rightarrow$ [B] $\Rightarrow$ [C] $\Rightarrow$ [D].
\end{lemma}

\noindent {\bf Proof}:

Clearly [A] $\Rightarrow$ [B].  Setting
$\mu_n = \delta(d_n)$ in [B] shows that [B] $\Rightarrow$ [C].

To show
that [C] $\Rightarrow$ [D], suppose that [C] holds but that
[D] fails to hold for some
compact $K \subset X$ and some $\epsilon > 0$.
Then there exists an $\epsilon > 0$, two increasing sequences of natural
numbers
$\{ n(k) \}$ and $\{ m(k) \}$ with
\[ n(k+1) > m(k) > n(k) \]
for all $k$,  and a sequence of points
$d_k \in D_{n(k)} \cap K$, such that
\begin{equation}
   \left| f_{m(k)}(d_k) - f_{n(k)}(d_k)\right| \ge \epsilon .
\label{contradict me}
\end{equation}
Since $K$ is compact, there exists an increasing sequence of natural numbers
$\{ k(j) \}$ such that  $\{ d_k(j) \}_{j=1}^\infty$ converges.
Now define the convergent sequence $\{ e_i \in D_{n(k(i))} \}_{i=1}^\infty$ by
$e_i = d_k(j)$ when $n(k(j)) \le i < n(k(j+1))$.  By [C], the
sequence $\{ f_i(e_i) \}$ converges.  But
$\{ f_i(e_i) \}$ {\it does not converge} along the subsequence indexed by
\[
    n(k(1)),\quad  m(k(1)), \quad n(k(2)),\quad  m(k(2)), \ldots
\]
 because of (\ref{contradict me}) and the fact that
\[
     e_{n(k(j))} = e_{m(k(j))} = d_{k(j)}.
\]
  This contradiction shows that [C] must imply [D].
\qquad $\blacksquare$

\begin{lemma}
\label{lemma2}
If condition [C] holds, then whenever $\{d_n \in D_n \}$ converges to
$x \in X$, the limit of $\{ f_n(d_n) \}$ depends only on $x$.
The function of $x$ which may thus be defined as
\[
 f(x) := \lim_{n \rightarrow \infty} f_n(d_n)
\]
when $d_n \longrightarrow x$, is continuous.
\end{lemma}

\noindent {\bf Proof}:

Assume that condition [C] holds.

Suppose $\{d_n \in D_n \}$ and $\{e_n \in D_n \}$ are two sequences
that both converge to $x \in X$.  Then the sequence
$d_1,e_2,d_3,e_4,d_5,e_6,\ldots$
also converges to $x$, and its $n^{th}$ term
is a member of $D_n$.  By condition [C], the sequence
$f_1(d_1), f_2(e_2), f_3(d_3), \ldots $ converges.  This shows that
$\lim f_n(d_n) = \lim f_n(e_n)$.

Suppose $x_k \longrightarrow x$ in $X$.  Given $\epsilon > 0$,
it is possible to find an increasing sequence $\{ n(k) \}$ of
natural numbers and a sequence of points $\{ e_{n(k)} \in D_{n(k)} \}$
such that $d_X(e_{n(k)},x_k) < {1 \over k}$ while
$  |f_{n(k)}(e_{n(k)}) - f(x_k)| < \epsilon $.
Then $\{ e_{n(k)} \}$ converges to $x$ just as $\{ x_k \}$ does, so
$\lim_{k \rightarrow \infty} f_{n(k)} (e_{n(k)}) = f(x)$.
Now
\[ |f(x_k) - f(x)| \le |f(x_k) - f_{n(k)}(e_{n(k)})|
    +   |f_{n(k)}(e_{n(k)}) - f(x)|  . \]
Since the last term tends to zero,
\[  \limsup_{k \rightarrow \infty} |f(x_k) - f(x)| \le \epsilon.  \]
Since $\epsilon$ may be arbitrarily small, $f(x_k) \longrightarrow f(x)$,
which shows that $f$ is continuous.
\qquad $\blacksquare$

\begin{lemma}
\label{lemma3}
If $(X,d_X)$ is a complete and separable metric space, and the functions
$f_n$ are bounded uniformly in $n$, then conditions
[A], [B], and [C] are all equivalent.
\end{lemma}
\noindent {\bf Proof}:

It remains to show that [C] $\Rightarrow$ [A] when $(X,d_X)$ is complete
and separable, and $\sup_{d \in D_n}\{ |f_n(d)| \} \le B$ for all $n$.

 Suppose that $\{ \mu_n \}$ converges to $\mu \in \P(X)$, where
$\mu_n(X \setminus D_n) = 0 $.  Since $X$ is complete and
separable, Prohorov's theorem implies that $\{ \mu_n \}$ is tight.
Thus, given $\epsilon > 0 $, there exists  a compact
$K_{\epsilon} \subset X$ such that $\mu_n (X \setminus
K_{\epsilon}) < \epsilon$ for all $n$.
With $f:X \longrightarrow {\mathbb R}$ as defined in Lemma \ref{lemma2},
\begin{eqnarray*}
   \left|\int_X f_n(x) \mu_n(dx) - \int_{X} f(x) \mu (dx) \right|
   \le && \\
    \left|\int_{K_{\epsilon}} f_n(x) - f(x) \mu_n(dx)\right| &+&
    2B\epsilon \\
     &+& \left|\int_{X} f(x) \mu_n(dx) - \int_{X} f(x) \mu (dx)\right|   .\\
\end{eqnarray*}
Condition [C] implies condition [D], a sort of uniform convergence
on compact sets that entails that
\[ \lim_{n \rightarrow \infty}
    \int_{K_{\epsilon}} \left|f_n(x)-f(x)\right| \mu_n(dx) = 0.   \]
Therefore,
\[ \limsup_{n \rightarrow \infty}
   \left|\int_{X} f_n(x) \mu_n(dx)
   - \int_{X} f(x) \mu (dx)\right| \le 2B\epsilon.
\]
Since $\epsilon$ is arbitrarily small, it follows that
\[    \int_{X} f_n(x) \mu_n(dx)
\longrightarrow \int_{X} f(x) \mu (dx).  \qquad \blacksquare  \]

\section{Theorems}
\label{Theorems}

Let $S$ and $T$ be separable metric spaces.  For each natural
number $n$,
let $K_n$ be a Markov transition from $S^n$ to $T^n$ that satisfies
the permutation condition (\ref {permute}).  Let $H_n$ be the transition
from $\e_n(S^n)$ to $\P(T)$ that is induced by $K_n$, as defined in
(\ref{InducedTransition}).

\begin{theorem}
\label{theorem1}
If a sequence of Markov transitions $\{ K_n \}$
propagates chaos, then there exists a continuous function
\[
      F:\P(S) \longrightarrow \P(T)
\]
such that,
if $ \e_n(\bs_n) \longrightarrow p $ in $\P(S)$ with $\bs_n \in S^n$, then
\[
     \left\{ \Kt_n(\bs_n,\cdot) \right\}_{n=1}^{\infty}
\]
is $F(p)$-chaotic.
\end{theorem}

\noindent {\bf Proof}:

The arguments of Lemma \ref{lemma1} and Lemma
\ref{lemma2} will be adapted to prove this.

Take $\P(S)$ with one of the metrics for the weak topology
to be the metric space $(X,d_X)$ of those lemmas, and take $\e_n(S^n)$
to be $D_n$.  For each bounded and continuous function $\phi \in C_b(\P(T))$
define the functions $ \phat_n: \e_n(S^n) \longrightarrow {\mathbb R}$ by
\begin{equation}
\label{phat}
     \phat_n(\zeta) := \int_{\P(T)} \phi(\eta) H_n(\zeta,d\eta) ,
\end{equation}
 where $H_n$ is as defined in (\ref{InducedTransition}).
These functions $\phi_n$ will play the role of the functions $f_n$ of
the lemmas.

By hypothesis, $\{K_n\}$ propagates chaos.
Proposition \ref{PoCDef2} therefore implies that
whenever $\{ \mu_n \in \P(\e_n(S^n)) \}_{n=1}^\infty$ converges in
$\P(\P(S))$ to
$\delta(p)$, then
\begin{equation}
\label{thingy}
   \int_{\e_n(S^n)} \phat_n(\zeta) \mu_n(d\zeta) \longrightarrow \phi(q)
\end{equation}
for some $q \in \P(T)$.  In fact, Proposition \ref{PoCDef2} implies
that $q$ does not depend on our choice of $\phi$: the same $q$ works
for all $\phi$ in (\ref{thingy}).

Condition (\ref{thingy}) resembles condition [B] of Lemma \ref{lemma1}.  Lemma
\ref{lemma1} and Lemma \ref{lemma2} can now be applied to show that there
exists a
continuous function $G_{\phi}(p)$, depending on $\phi$, such that if
 $\{ \bs_n \in S^n \}$ is a sequence satisfying
$ \e_n(\bs_n) \longrightarrow p $ in $\P(S)$, then
\[
     \phat_n(\e_n(\bs_n))  \longrightarrow G_{\phi}(p).
\]
By (\ref{thingy}), $G_{\phi}(p) = \phi(q)$ for some $q \in \P(T)$ {\it
that does not depend on} $\phi$.  The only way that all the
$G_{\phi}$'s can have this form and yet all be continuous is for the
dependence of $q$ on $p$ to be continuous: there
must be a continuous $F$ from $\P(S)$ to $\P(T)$ such that
$G_{\phi}(p) = \phi(F(p))$ for all $\phi \in C_b(\P(T))$.

Thus, there exists a continuous function $F$ from $\P(S)$ to $\P(T)$
such that
\[
    \left[ \e_n(\bs_n) \longrightarrow p \right] \implies
    \left[ \phat_n(\e_n(\bs_n))  \longrightarrow \phi(F(p)) \right]
\]
for all $\phi \in C_b(\P(T))$.  This fact, and the definitions
(\ref{InducedTransition}) and (\ref{phat}) of $H_n$ and $\phat$, imply
that
\[
    \left[ \e_n(\bs_n) \longrightarrow p \right] \implies
    \left[ \Kt_n(\bs_n, \cdot) \circ \e_n^{-1} \longrightarrow \delta(F(p))
\right].
\]
Finally, by Theorem \ref{SznitmanTanaka}, we have that
\[
    \left[ \e_n(\bs_n) \longrightarrow p \right] \implies
    \left\{ \Kt_n(\bs_n, \cdot)\right\}_{n=1}^{\infty} \  \mathrm{is} \
     F(p)\mathrm{-chaotic} . \qquad \blacksquare
\]

When $(S,d_S)$ is complete and separable, the necessary condition
of Theorem \ref{theorem1} is also {\it sufficient}.

\begin{theorem}[Main Theorem]
\label{theorem2}
Suppose $(S,d_S)$ is a complete, separable metric space.  Then $\{ K_n \}$
propagates chaos if and only if there exists a continuous function
\[ F:\P(S) \longrightarrow \P(T)  \]
such that, whenever $ \e_n(\bs_n) \longrightarrow p $ in $\P(S)$ with
$\bs_n \in S^n$, then
\[
     \left\{ \Kt_n(\bs_n,\cdot) \right\}_{n=1}^{\infty}
\]
is $F(p)$-chaotic.

\end{theorem}

\noindent {\bf Proof}:

We have just demonstrated that the condition is necessary (Theorem
\ref{theorem1}).
Next we demonstrate its sufficiency:

Suppose $P_n \in \P(S^n)$ is $p$-chaotic.  Let $\mu_n = P_n \circ \e_n^{-1}$.
Then $\mu_n \longrightarrow \delta(p)$ in $\P(\P(S))$ by Theorem
\ref{SznitmanTanaka}.   Our goal is to prove that
\[
    \int_{\P(S)} H_n(\zeta, d\eta) \mu_n(d\zeta) \longrightarrow
    \delta(F(p)),
\]
where $H_n$ is as defined in (\ref{InducedTransition}).
This is enough, by Proposition \ref{PoCDef2}, to demonstrate that
chaos propagates.

By hypothesis, if $\{\bs_n \in S^n\}$ is such that $\e_n(\bs_n)$
converges to $p$ then $\{ \Kt_n (\bs_n,\cdot) \}$ is $F(p)$-chaotic.
By Theorem \ref{SznitmanTanaka} and the fact that
\[
    \Kt_n(\bs_n, \cdot)\circ \e_n^{-1} = H_n(\e_n(\bs_n),\cdot) ,
\]
the hypothesis is equivalent to the statement that, if
$ p_n \in \e_n(S^n)$ for each $n$, then
\begin{equation}
\label{hypo2}
     \left[ p_n \longrightarrow p \right]
     \implies \left[H_n(p_n,\cdot) \longrightarrow \delta(F(p))
     \right].
\end{equation}

Let $\phi \in C_b(\P(T))$
be a bounded and continuous function on $\P(T)$, and define functions
$ \phat_n: \e_n(S^n) \longrightarrow {\mathbb R}$ by
\begin{equation}
\label{eq1}
     \phat_n(\zeta) := \int_{\P(T)} \phi(\eta) H_n(\zeta,d\eta) .
\end{equation}
These functions are uniformly bounded in $n$ since $\phi$ is bounded.

The hypothesis (\ref{hypo2}) and equation (\ref{eq1}) imply that
\begin{equation}
\label{eq2}
   \lim_{n \rightarrow \infty} \phat_n(p_n) =
   \phi(F(p))
\end{equation}
when $p_n \longrightarrow p$ with $ p_n \in \e_n(S^n)$.
We are assuming $S$ is complete and separable, therefore so is
$\P(S)$ \cite{Billingsley}.

We may now apply Lemma \ref{lemma3} with
\[
\begin{array}{ccc}
      X = \P(S),    &   D_n = \e_n(S^n), &  f_n = \phat_n,   \\
\end{array}
\]
and conclude that
\begin{equation}
\label{eq3}
  \lim_{n \rightarrow \infty} \int_{\P(S)} \phat_n(\zeta) \mu_n(d\zeta) =
          \phi(F(p))
\end{equation}
for any sequence $\{ \mu_n  \}$ that converges to
$\delta(p)$ in $\P(\P(S))$.

By equations (\ref{eq3}) and (\ref{eq1}),
\begin{eqnarray*}
 \phi(F(p)) &=&
 \lim_{n \rightarrow \infty} \int_{\P(S)} \phat_n(\zeta) \mu_n(d\zeta)  \\
 &=&
 \lim_{n\rightarrow \infty}\int_{\P(S)}\int_{\P(T)} \phi(\eta) H_n(\zeta,d\eta)
		    \mu_n(d\zeta)                      \\
 &=& \lim_{n \rightarrow \infty} \int_{\P(T)} \phi(\eta)
	    \int_{\P(S)} H_n(\zeta, d\eta) \mu_n(d\zeta)    ,  \\
\end{eqnarray*}
for all $\phi \in C_b(\P(T))$.
This implies that
\[
    \int_{\P(S)} H_n(\zeta, d\eta) \mu_n(d\zeta) \longrightarrow
    \delta(F(p))
\]
in $\P(\P(S))$, completing the proof.
\qquad $\blacksquare$

Theorem \ref{theorem1} states that the
{\it limit-law map} $F:\P(S) \longrightarrow \P(T)$
must be continuous.  If $S$ is complete and separable then,
conversely, any continuous map
$F$ is a possible limit-law map.  This fact is a corollary of Theorem
\ref{theorem2}:

\begin{corollary}

Suppose $(S,d_S)$ is a complete, separable metric space.  Then,
for any continuous $F:\P(S) \longrightarrow \P(T)$, there exists
a sequence of Markov transitions $\{ K_n \}_{n=1}^\infty$ that propagates
chaos,
and for which
\[
   \left\{ K_n(\bs_n,\cdot) \right\}_{n=1}^{\infty}
\]
is $F(p)$-chaotic whenever $ \e_n(\bs_n) \longrightarrow p $ in $\P(S)$.
\end{corollary}

\noindent {\bf Proof}:

Let $F:\P(S) \longrightarrow \P(T)$ be continuous.  For each $n$ and each
${\bs} \in S_n$, let
$K_n({\bs},\cdot)$ be the $n$-fold product measure
\begin{equation}
    K_n({\bs},\cdot) :=  F(\e_n({\bs}))\otimes F(\e_n({\bs})) \otimes
	  \cdots \otimes F(e_n({\bs})).
\end{equation}

Suppose the points $\bs_n \in S^n$ are such that $\e_n(\bs_n)$ converges
to $p$ as $n$ tends to infinity.
Since $F$ is continuous ,
$F(\e_n(\bs_n))$ converges to $F(p)$ as well,
so it is clear from Definition \ref{chaos} that the sequence of symmetric
measures $\left\{ K_n(\bs_n,\cdot)\right\}$ is
$F(p)$-chaotic.
By Theorem \ref{theorem2}, $\{ K_n \}$ propagates chaos.
\qquad $\blacksquare$

The Markov transition functions $\left\{ K_n \right\}$ may well be
{\it deterministic}, that is, the $n$-particle dynamics may simply be
given by a point-transformation from $S^n$ to $T^n$.   These
point-transformations are measurable maps from $S^n$ to $T^n$ that
commute with permutations of coordinates.

Let $f_n:S^n \longrightarrow T^n$ be a measurable map
that commutes with permutations of $n$-coordinates, i.e., such that
\begin{equation}
\label{permute2}
       f_n(s_{\pi(1)},s_{\pi(2)},\ldots,s_{\pi(n)}) = \pi \cdot
               f_n(s_1,s_2,\ldots,s_n)
\end{equation}
for each point ${\bf s} \in S^n$ and each
permutation $\pi$ of the symbols ${1,2,\ldots,n}$.
Given $f_n$, define the Markov transition $K_n$ from $S^n$ to $T^n$ by
\[
       K_n({\bf s},E) = {\bf 1}_E(f_n({\bf s}))
\]
when ${\bf s} \in S^n$ and $E \in \mathcal{B}_{T^n}$.
Say that $\{ f_n \}_{n=1}^{\infty}$ {\it propagates chaos} if the
sequence of deterministic
transition functions $\{ K_n \}$ propagates chaos.

  The following is an immediate corollary of Theorem \ref{theorem2}:

\begin{corollary}[Deterministic Case]
\label{deterministic}
Let $S$ be a Polish space, and for each $n$ let
$f_n$ be a measurable map from $S^n$ to $T^n$ that commutes with
permutations as in (\ref{permute2}).

$\left\{ f_n \right\} $ propagates chaos if and only if there exists a
continuous function
\[
      F:\P(S) \longrightarrow \P(T)
\]
such that $  \e_n(f_n({\bf s}_n)) \longrightarrow F(p) $
in $\P(T)$ whenever $ \e_n({\bf s}_n) \longrightarrow p $
in $\P(S)$.
\end{corollary}

\chapter{Conclusion}

We have studied the propagation of chaos by families of Markov
processes, having adopted a simple definition of propagation
of chaos, namely, that the processes propagate {\it all} chaotic sequences
of initial laws
to chaotic sequences.  Authors who wish to prove that certain families
of processes propagate chaos often show only that sequences of initial
laws of the form $\rho ^{\otimes n}$ are propagated to chaotic
sequences, that is, they show that pure chaos is propagated to chaos.
We have remarked that this does not imply unqualified propagation of
chaos.

Propagation of chaos, in its unqualified sense, entails the
continuity of the limit dynamics.  Families of Markov processes on
Polish spaces propagate chaos if and only if the associated Markov
transition functions satisfy the condition of Theorem \ref{theorem2}.

Our definition of propagation of chaos may be too simplistic to cover
some situations of interest.  For instance, the
subtle propagation of chaos phenomenon that is operative
in Lanford's validation of Boltzmann's equation --- where the chaos of
the initial laws propagates if those laws have densities that
converge uniformly --- is not subject to our treatment here.

Further foundational research on the
propagation of chaos phenomenon of Lanford's theorem is called for.
Is the phenomenon endemic to the Boltzmann-Grad limit, or is it, like
the propagation of chaos that is the subject of our theorems,
a more general probabilistic phenomenon that should appear in other parts
of kinetic theory?


\begin{thebibliography}{99}

\bibitem[1]{Billingsley} P. Billingsley. {\it Convergence of
Probability Measures}. John Wiley \& Sons, New York, 1968.

\bibitem[2]{Bird} G. A. Bird.  {\it Molecular Gas Dynamics}.
Clarendon Press, Oxford, 1976.

\bibitem[3]{Boltzmann}L. Boltzmann. {\it Lectures on Gas
Theory}. Dover Publications, New York, 1995.

\bibitem[4]{BossyTalay} M. Bossy and D. Talay.  A stochastic
particle method for the McKean-Vlasov and the Burgers equation. {\it
Mathematics of Computation} 66 (217): 157-192, 1997.

\bibitem[5]{BH}W. Braun and K. Hepp. The Vlasov dynamics and its
fluctuations in the $\oon$ limit of interacting classical
particles.  {\it Communications in Mathematical Physics} 56: 101-113,
1977.

\bibitem[6]{Chorin}A. J. Chorin.  Numerical study of slightly viscous
flow.  {\it Journal of Fluid Mechanics} 57: 785-796, 1973.

\bibitem[7]{CIP}C. Cercignani, R. Illner, and M. Pulvirenti.  {\it The
Mathematical Theory of Dilute Gases}.  Springer-Verlag, New York, 1994.

\bibitem[8]{DawsonGartner}D. Dawson and J. G\"artner.  Large deviations
from the McKean-Vlasov limit for weakly interacting diffusions.
{\it Stochastics} 20: 247-308, 1987.

\bibitem[9]{DipLio} R. Di Perna and P.L. Lions.  On the Cauchy problem
for the Boltzmann equation.  {\it Annals of Mathematics}
130: 321-366, 1989.

\bibitem[10]{Dudley}R. M. Dudley.  {\it Real Analysis and Probability}.
Wadsworth \& Brooks/Cole, Pacific Grove, California, 1989.

\bibitem[11]{DupuisEllis} P. Dupuis and R. Ellis.  {\it A Weak
Convergence Approach to the Theory of Large Deviations}.    John Wiley
\& Sons, New York, 1997.

\bibitem[12]{EthKur}S.N. Ethier and T.G. Kurtz. {\it Markov Processes:
Characterization and Convergence}.  John Wiley \& Sons, New York, 1986.

\bibitem[13]{Grad} H. Grad.  On the kinetic theory of rarefied gases.
{\it Communications in Pure and Applied Mathematics}   2: 331-407,
1949.

\bibitem[14]{Graham} C. Graham.  McKean-Vlasov Ito-Skorohod equations,
and nonlinear diffusions with discrete jump sets.  {\it Stochastic
Processes and their Applications}  40: 69-82, 1992.

\bibitem[15]{Gru}F. A. Gr\"unbaum.  Propagation of chaos for the
Boltzmann equation. {\it Archive for Rational Mechanics and
Analysis}  42: 323-345, 1971.

\bibitem[16]{Kac}M. Kac. {\it Probability and Related Topics in
Physical Sciences}.  American Mathematical Society, Providence, Rhode
Island, 1976.

\bibitem[17]{Kac55}M. Kac. Foundations of kinetic theory. {\it
Proceedings of the Third Berkeley Symposium on Mathematical
Statistics and Probability, Vol III}.  University of California Press,
Berkeley, California, 1956.

\bibitem[18]{Lanford} O.E. Lanford III.  The evolution of large
classical systems.  {\it Lecture Notes in Physics}, 35: 1-111.
Springer-Verlag, Berlin, 1975.

\bibitem[19]{Le Cam}L. Le Cam. {\it Convergence in Distribution of
Stochastic Processes}.  University of California Publications in
Statistics. 2 (11): 207-236, 1957.

\bibitem[20]{MarPul}C. Marchioro and M. Pulvirenti.  Hydrodynamics in
two dimensions and vortex theory.  {\it Communications in Mathematical
Physics}   84: 483-503, 1982.

\bibitem[21]{McK66}H. P. McKean, Jr.  A class of Markov processes
associated with nonlinear parabolic equations.  {\it Proceedings of
the National Academy of Science} 56: 1907-1911, 1966.

\bibitem[22]{McK68}H. P. McKean, Jr.  Propagation of chaos for a class
of nonlinear parabolic equations.  {\it Lecture Series in Differential
Equations}  7: 41-57.  Catholic University, Washington, D.C., 1967.

\bibitem[23]{McK75}H. P. McKean, Jr.  Fluctuations in the kinetic
theory of gases.  {\it Communications in Pure and Applied
Mathematics}  28: 435-455, 1975.

\bibitem[24]{Meleard}S. M\'el\'eard. Asymptotic behavior of some
interacting particle systems; McKean-Vlasov and Boltzmann models.
{\it Lecture Notes in Mathematics}, 1627. Springer-Verlag, Berlin, 1995.

\bibitem[25]{Meleard2} S. M\'el\'eard.  A probabilistic proof of the
vortex method for the 2D Navier-Stokes equations.  Preprint.

\bibitem[26]{Oelsch}K. Oelschl\"ager.  A law of large numbers for
moderately interacting diffusion processes.  {\it
Zeitschrift f\"ur Wahrscheinlichkeitstheorie und verwandte Gebiete}
69: 279-322, 1985.

\bibitem[27]{Osada}H. Osada.  Propagation of chaos for the two
dimensional Navier-Stokes equation. {\it Probabilistic Methods in
Mathematical Physics}.  Academic Press, Boston, 1987.

\bibitem[28]{Prohorov}Yu. V. Prohorov.  Convergence of random processes
and limit theorems in probability theory. {\it Theory of Probability
and it Applications} 1: 157-214, 1956.

\bibitem[29]{Pulvirenti} M. Pulvirenti. Kinetic limits of stochastic
particle systems.
{\it Lecture Notes in Mathematics}, 1627. Springer-Verlag, Berlin, 1995.

\bibitem[30]{Rez}F. Rezakhanlou.  Kinetic limits for a class of
interacting of interacting particle systems. {\it  Probability Theory and
Related
Fields} 104: 97-146, 1996.

\bibitem[31]{RezTar}F. Rezakhanlou and J. Tarver.  Boltzmann-Grad limit
for a particle system in continuum.  {\it Annales de l'Institut Henri
Poincar\'e Probab. Stat.}  33 (6): 753-796, 1997.

\bibitem[32]{ShigaTanaka}T. Shiga and H. Tanaka. Central limit theorem
for a Markovian system of particles in mean field interaction. {\it
Zeitschrift f\"ur Wahrscheinlichkeitstheorie und verwandte Gebiete}
69: 439-459, 1985.

\bibitem[33]{Spohn}H. Spohn. {\it Large Scale Dynamics of Interacting
Particles}. Springer-Verlag, Berlin, 1991.

\bibitem[34]{Szn}A. Sznitman. \'Equations de type de Boltzmann,
spatialement homog\`enes.  {\it Zeitschrift f\"ur
Wahrscheinlichkeitstheorie und verwandte Gebiete} 66: 559-592, 1984.

\bibitem[35]{Sz84}A. Sznitman.  Nonlinear reflecting diffusion
process, and the propagation of chaos and fluctuations associated.
{\it Journal of Functional Analysis} 56: 311-336, 1984.

\bibitem[36]{Sznitman}A. Sznitman. Topics in propagation of chaos.
{\it Lecture Notes in Mathematics}, 1464. Springer-Verlag, Berlin,
1991.

\bibitem[37]{Talay}D. Talay.  Probabilistic numerical methods for
partial differential equations: elements of analysis.
{\it Lecture Notes in Mathematics}, 1627. Springer-Verlag, Berlin,
1995.

\bibitem[38]{Tanaka} H. Tanaka.  Limit theorems for certain diffusion
processes with interaction.  {\it Taniguchi Symposium on Stochastic
Analysis} (pp. 469-488). Katata, Kyoto, 1982.

\bibitem[39]{Uch83}K. Uchiyama.  Fluctuations of Markovian systems in
Kac's caricature of a Maxwellian gas.  {\it Journal of the Mathematical
Society of Japan} 35 (3): 477-499, 1983.

\bibitem[40]{Uch88}K. Uchiyama. Fluctuations in a Markovian system of
pairwise interacting particles. {\it Probability Theory and Related
Fields} 79: 289-302, 1988.



\end{thebibliography}
\end{document}